\documentclass[twocolumn]{svjour3}
\usepackage[utf8]{inputenc}
\usepackage{epsfig} 
\usepackage{amsfonts,latexsym,amssymb} 
\usepackage{mathrsfs,amsmath}
\usepackage{graphicx}
\usepackage{xcolor}
\usepackage{hyperref}
\usepackage{caption}
\usepackage{subcaption}
\usepackage{tikz}

\usepackage{mathptmx}
\usetikzlibrary{arrows,matrix}
\usepackage{epstopdf}
\usepackage{multirow}
\usepackage{array}
\usepackage{hyperref}
\usepackage{amsmath}
\usepackage{lineno}
\usepackage{cite}

\DeclareMathAlphabet      {\mathbfit}{OML}{cmm}{b}{it}

\newcommand{\vc}{\mbox{${\mathbf{c}}$}}

\def \begpf{ \begin{description} \item[{\bf Proof:\hspace{1.5em} }]  }
\def \endpf{ \begin{flushright} $\Box $ \end{flushright}   \end{description}}

    \newtheorem{prop}{Proposition}[section]
    
    \newtheorem{lemm}[prop]{Lemma}

    \newtheorem{exmpl}{Example}
    
    \newcommand{\begex}[2]{ \begin{exmpl}[#1] \label{#2} \end{exmpl}
    \begin{description}  \item
   }
    \def  \endex {       \end{description}  }

\def \begpf{ \begin{description} \item[{\bf Proof:\hspace{0.75em} }]  }
\def \endpf{ \begin{flushright} $\Box $ \end{flushright}   \end{description}}
\def \begrm{ \begin{description} \item[{\bf Remark:\hspace{0.5em} }]  }

\def \endrm{  \end{description} }
\def \begex{ \begin{description} \item[{\bf Example:\hspace{0.5em} }]  }
\def \endex{  \end{description} }

\newcommand{\beq}[1]{ \begin{equation} \label{eq.#1} }
\newcommand{\eeq}{ \end{equation} }

\newcommand{\blemm}[1]{ \begin{lemm} \label{lemm.#1} }
\newcommand{\elemm}{ \end{lemm} }
\newcommand{\barr}{ \begin{array}} 
\newcommand{\earr}{ \end{array} }

\begin{document}
\title{\centering Passive Vibration-Driven Locomotion}
\author{Anna Zigelman$^{1,\ast}$\thanks{$^{\star}$Corresponding Author \email{annar@technion.ac.il}}  \and Gilad Israel$^{1}$ \and Yizhar Or$^{1}$ \and Yuli Starosvetsky$^{1}$}


\institute{
	1. Faculty of  Mechanical Engineering,
	Technion - Israel Institute of  Technology, Haifa 3200003, Israel.
}

\maketitle

\begin{abstract}
We investigate a concept of passive, vibration-driven locomotion, in which a mechanical system achieves horizontal self-propulsion by resonantly harvesting energy from vertical environmental excitations (e.g. ambient vibrations of underwater pipelines), without a direct propulsive actuation. The system consists of a capsule containing an internal pendulum attached to its base mounted on a vertically vibrating substrate. The underlying locomotion mechanism relies on resonant energy transfer from the vertically vibrating substrate to the internal oscillatory element. Under appropriate forcing conditions and in the presence of asymmetric dissipative interactions, this internal oscillator induces a net unidirectional motion of the capsule. The analysis focuses on regimes of progressive motion arising in the vicinity of parametric resonances. Two asymptotic limits are considered: small-amplitude parametric excitation leading to a (2:1) resonant oscillatory motion of the pendulum, and large-amplitude excitation leading to a (1:1) resonant unidirectional rotational motion of the pendulum. Given the asymmetry of the dissipative force acting on the capsule, both resonant regimes result in a progressive motion of the capsule system. To identify optimal locomotion regimes in both cases, we employ tailored asymptotic approaches based on multi-scale expansions and direct averaging analysis. The resulting slow-flow and averaged-flow models reveal the full bifurcation structure of steady-state solutions associated with forward capsule motion for both low- and high- amplitude excitations. Analytical predictions are shown to be in good agreement with direct numerical simulations of the full \linebreak capsule-pendulum system.
\end{abstract}

\keywords{Passive Vibration-Driven Locomotion, Parametric Resonance, Resonant Energy Transfer, Strongly Nonlinear Regimes}

%
%
%


\section{Introduction}\label{Intro} 
Vibration-driven (VD), locomotion constitutes a distinct \linebreak class of self-propelled systems in which progressive motion arises from oscillatory excitation externally stimulated on internal oscillatory components embedded within the capsule, while the capsule itself interacts with a resistive surrounding medium. In such systems, propulsion does not rely on continuous external thrust or on the cyclic actuation of protruding elements such as legs, wheels, bristles, or tracks (see e.g., \cite{Raibert_1986,Hirose,Collins_2005,Holmes_2006,Raibert_2008,Hwangbo_2019,Bledt_2018,Tagliavini_2022,Bruzzone_2022,Supik_2023,Fath_2024}). Instead, the net locomotion emerges from the externally stimulated oscillations of the internal vibrating device along with direction-dependent resistive (dissipative) forces exerted by the surrounding medium on the capsule during the oscillation cycle. A canonical realization of this concept is the capsule-type architecture, typically modeled as one outer shell or multiple coupled outer shells enclosing one or more internal oscillatory components, whose motion alters the instantaneous velocity of the capsule and, consequently, the resistive forces acting on it. By appropriately shaping the internal dynamics and properly manipulating the external resistive forces, oscillatory motion can be converted into a sustained drift of the capsule through the medium \cite{Li_2006,CHERNOUSKO2008116,FANG20114002,Liu_2013,Zimmermann_2009,FANG2017153,Du_2018,Xu_2019}. Compared to locomotion strategies  based on protruding appendages, vibration-driven capsule systems offer several fundamental advantages. Their sealed and compact geometry enhances robustness and reliability, particularly in confined, cluttered, or hostile environments where protruding elements may suffer from entanglement, wear, or mechanical failure. Moreover, the absence of external appendages simplifies the mechanical design and facilitates miniaturization, which has motivated extensive research on capsule-based locomotion in medical, industrial, and \linebreak inspection-related applications~\cite{Carpi_2007,Ciuti_2011,Nelson_2010,Zhang_2011,ALHASSAN2020105862}.

Closely related physical principles also underlie a complementary class of vibration-induced transport mechanisms in which the oscillatory excitation is supplied externally \linebreak rather than generated internally. In this context, the term passive vibration-driven locomotion refers to systems in which net translational motion emerges directly from externally imposed mechanical vibrations of the environment, without internal actuation or active control (see e.g. \cite{Nath_2022,Viswarupachari_2012,Barois_2024,Dorbolo_2005,Xu_2017} and references therein). A canonical realization of this concept is provided by bristlebot-type devices placed on vibrating tables, where directional drift arises from asymmetric frictional contacts, compliance of the bristles, and intermittent stick-slip dynamics under zero-mean excitation \cite{Supik_2023,Fath_2024,Cicconofri_2016,Noselli_2014}. Related mechanisms have subsequently been explored in a broader class of passive vibration-driven systems. These include the rigid or compliant bodies interacting with vibrating substrates, where nonlinear contact dynamics, geometric asymmetry, and internal physical constraints similarly convert oscillatory excitation into directed translational or rotational motion~\cite{Hatatani_2025,Broseghini_2019,Baule_2013}. More recent studies have further extended these ideas to structurally constrained or guided configurations, in which passive bodies (i.e. the bodies containing no internal actuation) interact with externally vibrating substrates and confining geometries, and where base excitation, together with dry-friction dynamics, gives rise to directed transport along preferred directions~\cite{Kilikevičius_2021,El_Banna_2026,Kilikevičius_2022,Kilikeviius2022AnalysisOC,Kilikevičius_2021A}.

In the present work we introduce an additional concept of vibration-driven locomotion that lies at the interface between internally excited capsule locomotion and externally driven passive rectification. Specifically, we consider a passive vibration-driven capsule in which environmental or base vibrations do not directly translate the capsule, but instead resonantly excite an internal, fully passive \linebreak oscillatory/rotational mechanism. To the best of the authors’ knowledge, the use of a passive internal pendulum resonator for vibration-driven capsule locomotion has not been previously investigated. In addition, a comprehensive asymptotic analysis of both weakly and strongly nonlinear resonant regimes in capsule models with internally resonant parametric excitation has not been reported in prior studies. This internal mechanism is mechanically realized as a pendulum pivoted to the capsule and subject to vertical harmonic base excitation. The resulting parametrically excited resonant oscillations or rotations of the internal pendulum provide an internal “motion resource” that can be systematically rectified into progressive capsule motion through direction-dependent dissipation at the capsule-environment interface. To this end, we adopt a simple yet physically realizable model of asymmetric viscous damping, which enables progressive motion of the internally excited capsule. Such damping asymmetry can be implemented, for example, by introducing a controlled asymmetry in the capsule geometry, leading to different effective drag forces during forward and backward motion through the same dissipative medium (see, e.g., \cite{Benham_Boucher_Labbé_Benzaquen_Clanet_2019}). Similar models of asymmetric viscous resistance have been previously employed in foundational vibration-driven locomotion studies (see, e.g., \cite{Zimmermann_2010}). In parallel, the use of internal rotators or pendula connects the present capsule architecture with a broad literature on rotary nonlinear energy sinks, targeted energy transfer, and multi-dimensional (2D and 3D) energy channeling devices \cite{Vorotnikov_2018}. In particular, the concept of 2D and 3D energy channeling developed in \cite{Vorotnikov_2015,Vorotnikov_2015B,Vorotnikov_2018B,Jayaprakash_2017,Jayaprakash_2017B} explicitly exploits internal rotational elements to achieve resonant, multi-directional energy routing. Finally, although the primary focus of the present study is locomotion, pendulum- and rotator-based resonance mechanisms also play a central role in contemporary designs of pendulum-based vibration energy harvesters, underscoring the broader relevance of internal resonance and nonlinear rectification under broadband excitations \cite{WANG2023}. 

In the present study, we asymptotically analyze the emergence of distinct locomotion states of the capsule associated with two different internal resonance mechanisms: (i) a (2:1) internal parametric resonance manifested by the weakly nonlinear oscillations of the pendulum, and (ii) a (1:1) resonance giving rise to sustained rotational motion of the same internal pendulum. To elucidate the mechanisms underlying the formation of these locomotion states and to develop predictive analytical tools for their excitation and control, we examine the system dynamics at two complementary and robust asymptotic levels. Specifically, we analyze the weakly nonlinear dynamics of the pendulum in the vicinity of the (2:1) parametric resonance, as well as the strongly nonlinear rotational dynamics in the neighborhood of the (1:1) resonant regime. As demonstrated in the paper, both vibrational regimes induce progressive motion of the capsule, and the corresponding mean translational velocities can be accurately approximated using the derived reduced-order analytical models. On the basis of these asymptotic results, we identify parameter domains that - guarantee, efficient and controllable locomotion regimes. As shown below, the parametric domains may contain multiple steady-state solutions, including both trivial and nontrivial states. The bifurcation structure of the possible locomotion regimes is revealed through asymptotic analysis in the limits of both low-amplitude excitation and strong excitation. The analytical predictions are shown to be in good agreement with direct numerical simulations of the full capsule model.

The structure of the paper is as follows. In Section~\ref{S:sec 2}, we derive the nondimensional equations of motion governing the capsule dynamics. Section~\ref{S:3} presents numerical evidence for the coexistence of two distinct locomotion states of the capsule. These states are then analyzed in detail in Section~\ref{S:41}, devoted to oscillatory response and {Section}~\ref{S:42}, which examines rotatory response. Finally, Section~\ref{S:6} summarizes the main findings of the present study and outlines directions for future research.

\section{Model and formulation}\label{S:sec 2}
\begin{figure}[ht!]
\includegraphics[width=0.5\textwidth]{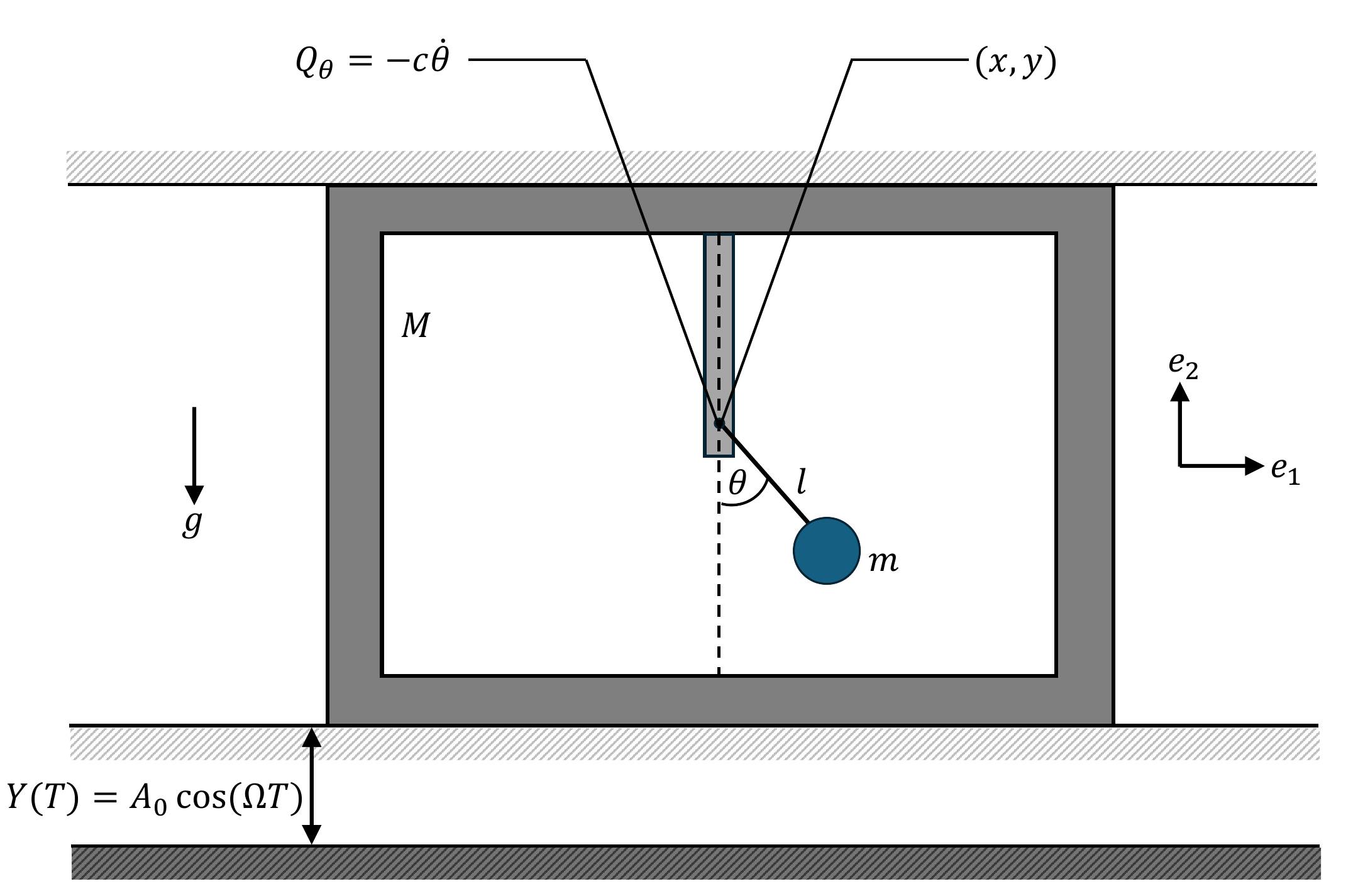}
\caption{Schematic of the passive capsule-pendulum - model.}\label{F:fig1}
\end{figure}

The system under consideration is a 2D model comprising a capsule, which consists of a rectangular frame of mass $M$ and an internal pendulum of mass $m$ and length $l$. The system is subjected to a harmonic lateral base excitation. Schematic of the model is presented in Figure~\ref{F:fig1}. The main mass is constrained to a vertically vibrating base, and can move on the horizontal plane only, with no out-of-plane possible movement. The  tangential friction forces between the capsule and the substrate are subject to an asymmetric viscous damping, where the damping coefficient of the forward motion is smaller than that of the backward motion. The inner pendulum is subjected to an internal dissipative torque due to friction assumed at its hinge. This dissipative torque is assumed to be linearly proportional to the angular velocity of the pendulum. We make here an additional physical assumption that the mass of the pendulum is small in comparison to the mass of the capsule, and the mass of the pendulum link is negligible. Gravity is taken into account as well. We assume that the capsule starts from rest while the inner pendulum starts from certain initial conditions, and may exhibit large amplitude excursions. 

We choose $X$, $Y$, and $\theta$ as the generalized coordinates of the system which denote horizontal and vertical position of the capsule, respectively, as well as the angle of internal pendulum. Lagrangian of the underlying undamped system subject to a base excitation is given by,
\begin{equation}\label{E:1}
\begin{aligned}
    L &=  
\frac{1}{2}M\dot{Y}{^2} + \frac{1}{2}M\dot{X}{^2} + \frac{1}{2}m\bigl[\dot{Y}^2 + \dot{X}^2 + {l^2}\dot{\theta}^2\\
&+ 2l\dot{\theta}\left( \dot{Y}\sin (\theta ) + \dot{X}\cos (\theta ) \right) \bigr]\\
&- \left( {MgY + mg\left( {Y + l(1 - \cos (\theta )} \right)} \right),
\end{aligned}
\end{equation}
where $M$ denotes the mass of the capsule, $m$ is the tip-mass, $l$ is the length of the pendulum, and $g$ is the gravity acceleration. Moreover, $T$ denotes time, and we use the following notation for the temporal derivatives, $\dot{Y}=dY/dT$, etc.

The capsule is constrained to oscillatory motion in the vertical direction, with constant amplitude $A_0$  and constant frequency $\Omega$, namely
\begin{equation}\label{E:2}
    Y(T) = {A_0}\cos (\Omega T).
\end{equation}
Moreover, the capsule is subjected to asymmetric viscous damping, where the damping coefficient is \linebreak direction-dependent and varies in magnitude with respect to the sign of capsule's velocity,
\begin{equation}\label{E:3}
    {F_d} =  - \lambda \left( \dot{X} \right)\dot{X},\quad \lambda \left( \dot{X} \right) = \left\{ \begin{array}{l}
{\lambda _1}, \,\, \dot{X} > 0\\
{\lambda _2},\,\, \dot{X} \le 0
\end{array} \right.,\,\,\,{\lambda _{1,2}} > 0,
\end{equation}
where $\lambda(\dot{X})$ is linear asymmetric viscous damping acting on the capsule.

The linear dissipation torque at the hinge of the internal pendulum is given by
\begin{equation}\label{E:4}
    {Q_\theta } =  - c\dot{\theta},
\end{equation}
where $c$ is a linear damping coefficient of the internal pendulum.

Using the Euler-Lagrange formalism and accounting for the external and internal dissipation, we arrive at the following equations of motion governing the dynamics of the capsule – pendulum model,
\begin{equation}\label{E:5}
    \begin{aligned}
&\left( {M + m} \right)\ddot{X} + ml\ddot{\theta}\cos (\theta ) - ml\dot{\theta} ^2\sin (\theta ) =  - \lambda (\dot{X})\dot{X},\\
&ml\left[\ddot{X}\cos (\theta ) + l\ddot{\theta} + \left(g - \Omega ^2{A_0}\cos (\Omega T)\right)\sin (\theta )\right] =  - c\dot{\theta}.
\end{aligned}
\end{equation}

To express the system in~\eqref{E:5} in dimensionless form, we perform the following time and coordinate scalings and introduce the nondimensional parameters:

\begin{equation}\label{E:6}
    \begin{aligned}
&x = \frac{X}{l},\quad t = {\Omega _s}T, \quad {\Omega _s} = \sqrt {\frac{g}{l}}, \quad \omega  = \frac{\Omega }{{{\Omega _s}}},\\
&\mu \left( {\dot x} \right) = \sqrt {\frac{l}{g}} \frac{{\lambda \left( {\dot x} \right)}}{{\left( {M + m} \right)}} = \,\left\{ \begin{array}{l}
{\mu _1},\,\,\,\,\dot x \ge 0\\
{\mu _2},\,\,\,\,\dot x < 0{\rm{ }}
 \end{array} \right., \\
 & \zeta  = \frac{c}{{m{l^2}}}\sqrt {\frac{l}{g}}, \quad \varepsilon  = \frac{m}{{\left( {M + m} \right)}}, \quad {A} = \frac{{{\Omega ^2}{A_0}}}{g}. 
\end{aligned}
\end{equation}

Upon nondimensionalizing the governing equations of motion, the system in~\eqref{E:5} takes the following form:
\begin{equation}\label{E:7}
    \begin{aligned}
&x'' + \varepsilon \cos (\theta ) \theta'  - \varepsilon {{\theta'}^2}\sin (\theta ) =  - \mu (x')x',\\
&x''\cos (\theta ) + \theta''  + \left( {1 - {A}\cos (\omega t )} \right)\sin (\theta ) =  - \zeta\theta', \\
&\mu \left( {x'} \right) = \left\{ \begin{array}{l}
{\mu _1},\,\,\,\,x' > 0,\\
{\mu _2},\,\,\,\,x' \le 0,
\end{array} \right.
\end{aligned}
\end{equation}
where we use the notation $x'=dx/dt$, etc.

Note that the parameter $\varepsilon$, which according to~\eqref{E:6} has the physical meaning of a mass ratio (i.e., the ratio between the mass of the internal pendulum and the total system mass), is introduced as a formal, non-dimensional asymptotic scaling parameter and is assumed to be small, $0 < \varepsilon \ll 1$.

In Sections~\ref{S:3}--\ref{S:42}, we conduct an extensive analytical and numerical study of the nonstationary dynamical regimes exhibited by~\eqref{E:7}, focusing on the parametric excitation in the vicinity of the (2:1) oscillatory resonance response (see Section~\ref{S:41}), which means that the actuation frequency is doubled relative to the system’s natural frequency, as well as the (1:1) rotatory resonance response (see Section~\ref{S:42}), which means that the actuation frequency is equal to the natural frequency. 
Special emphasis is placed on the analytical characterization of the various resonant response regimes that enable unidirectional propulsion of the capsule.

\begin{figure*}[ht!]
\includegraphics[width=\textwidth]{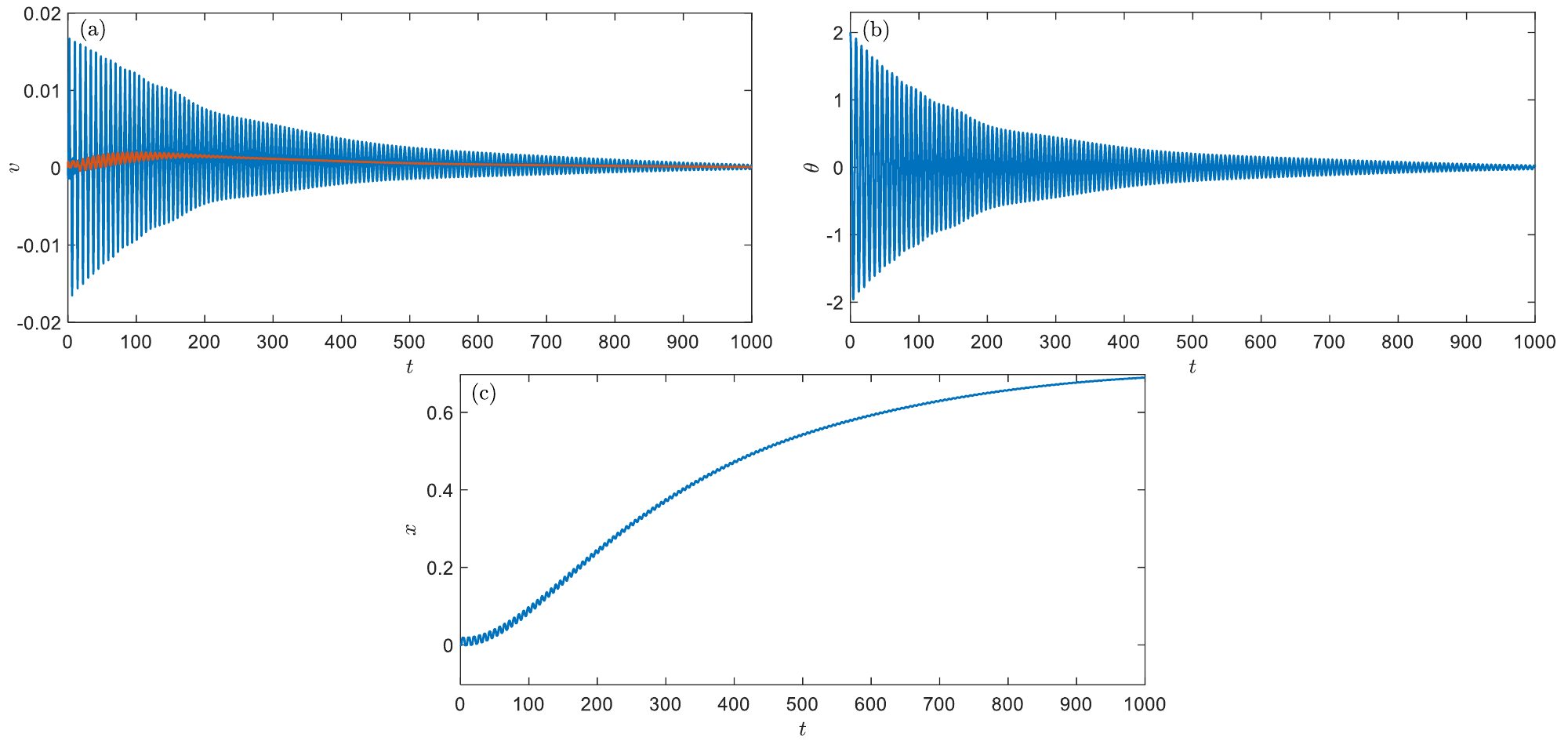}
\caption{(\textbf{Case 1}) Time histories illustrating the pendulum and capsule response under (2:1) resonance conditions, where the pendulum exhibits an oscillatory motion. System parameters: $\omega=2$, $\varepsilon=0.01$, $A=0.01$, $\zeta=0.01$, $\mu_1=0.01$, and $\mu_2=0.02$. Initial conditions: $x(0)=0$, $\theta(0)=2$, $x'(0)=0$, and $\theta'(0)=0$. (a) The time history of the capsule's velocity $v=x'$ versus time $t$, (b) the time history of the pendulum's angle $\theta$ versus time, and (c) the time histories of capsule's position $x$ versus time. 
}\label{F:fig11}
\end{figure*}

\begin{figure*}[ht!]
\includegraphics[width=\textwidth]{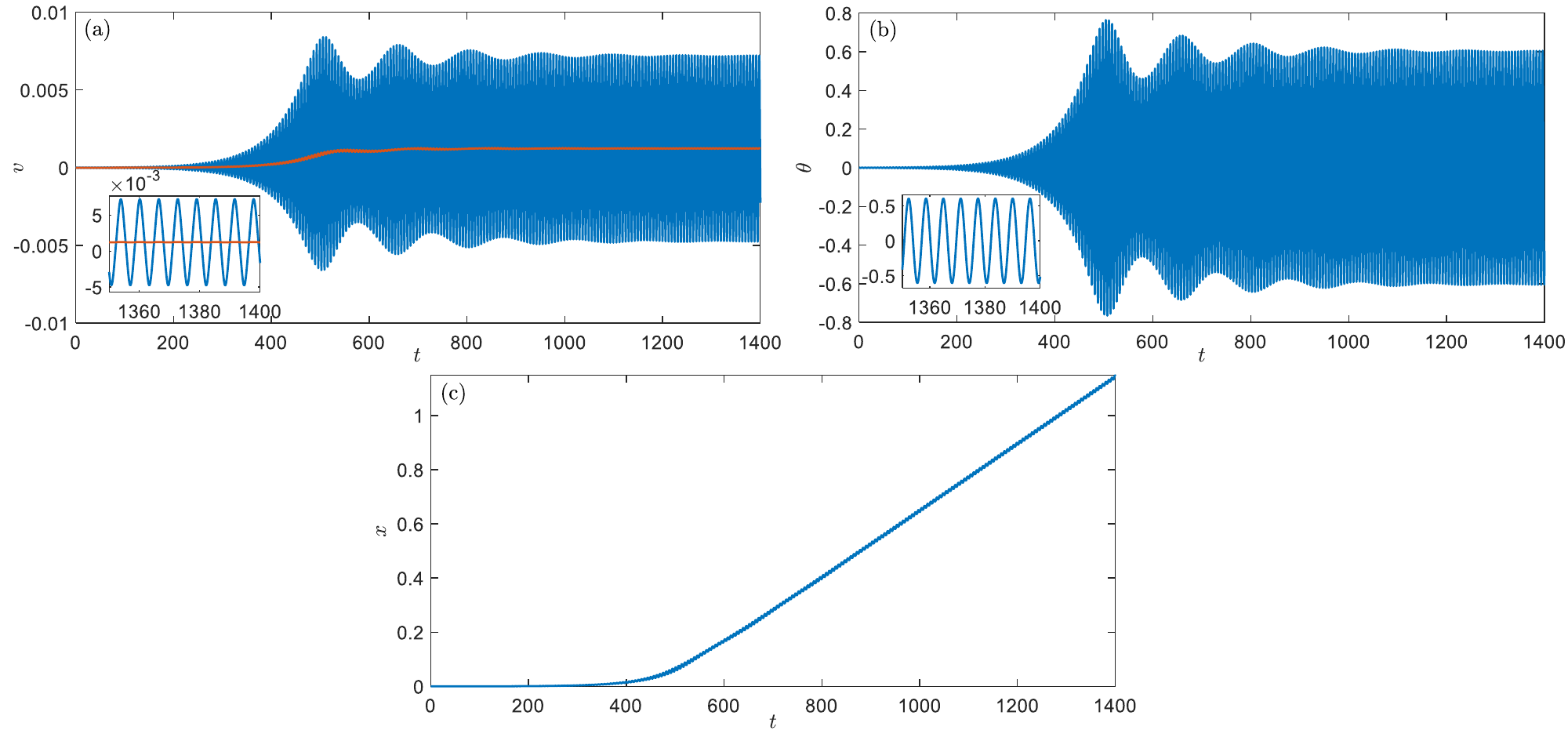}
\caption{(\textbf{Case 2}) Time histories illustrating the pendulum and capsule response under (2:1) resonance conditions, where the pendulum exhibits an oscillatory motion. System parameters: $\omega=2$, $\varepsilon=0.01$, $A=0.08$, $\zeta=0.01$, $\mu_1=0.01$, and $\mu_2=0.02$. Initial conditions: $x(0)=0$, $\theta(0)=0.001$, $x'(0)=0$, and $\theta'(0)=0$. (a) The time history of the capsule's velocity $v=x'$ versus time $t$, (b) the time history of the pendulum's angle $\theta$ versus time, and (c) the time histories of capsule's position $x$ versus time. The insets in panels (a) and (b) show the magnified views of the corresponding graphs, and the orange curve in panel (a) depicts the running average of the velocity.}\label{F:fig2}
\end{figure*}

\begin{figure*}[ht!]
\includegraphics[width=\textwidth]{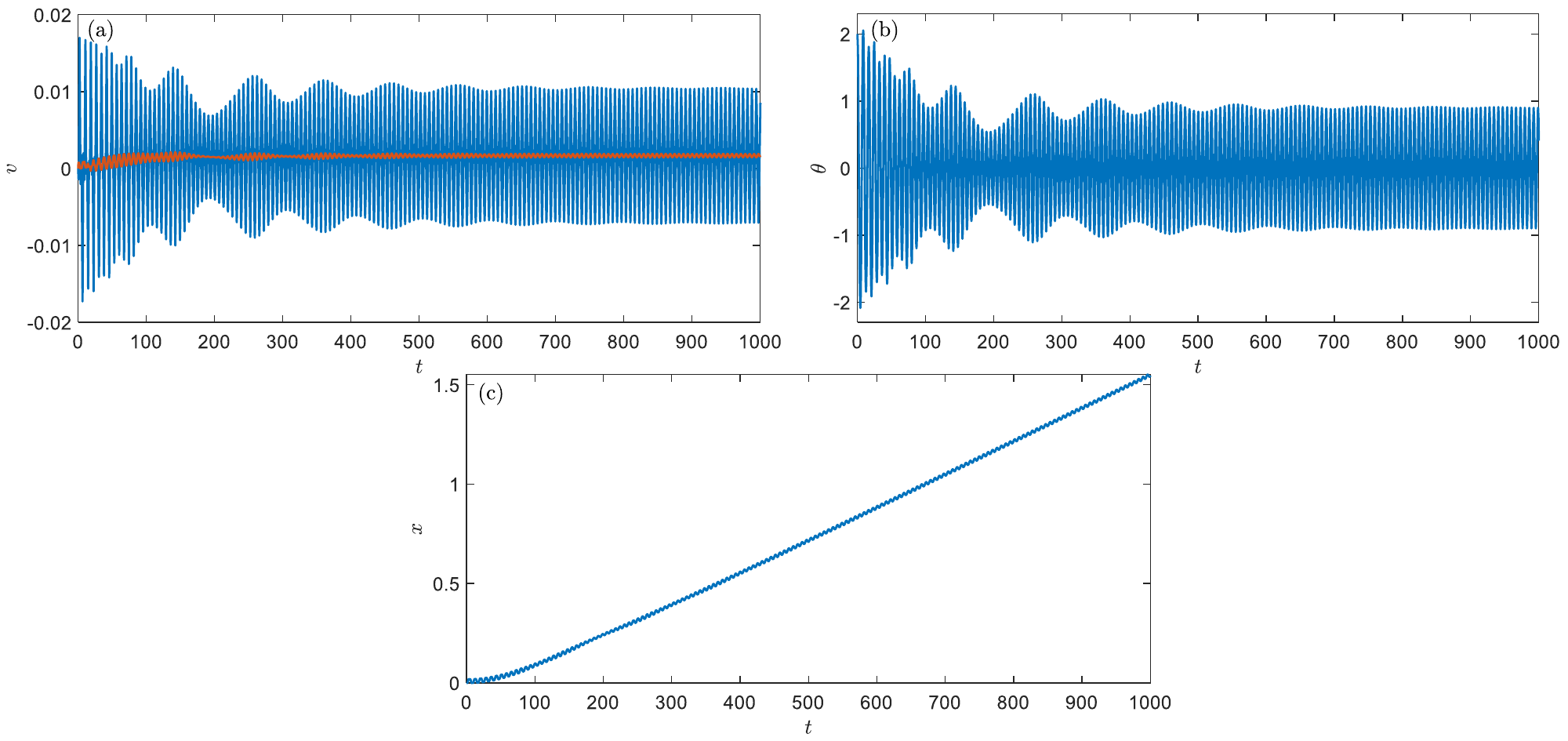}
\caption{(\textbf{Case 3}) Time histories illustrating the pendulum and capsule response under (2:1) resonance conditions, where the pendulum exhibits an oscillatory motion. System parameters: $\omega=1.94$, $\varepsilon=0.01$, $A=0.08$, $\zeta=0.01$, $\mu_1=0.01$, and $\mu_2=0.02$. Initial conditions: $x(0)=0$, $\theta(0)=2$, $x'(0)=0$, and $\theta'(0)=0$. (a) The time history of the capsule's velocity $v=x'$ versus time $t$, (b) the time history of the pendulum's angle $\theta$ versus time, and (c) the time histories of capsule's position $x$ versus time.}\label{F:fig12}
\end{figure*}

\begin{figure*}[ht!]
\includegraphics[width=\textwidth]{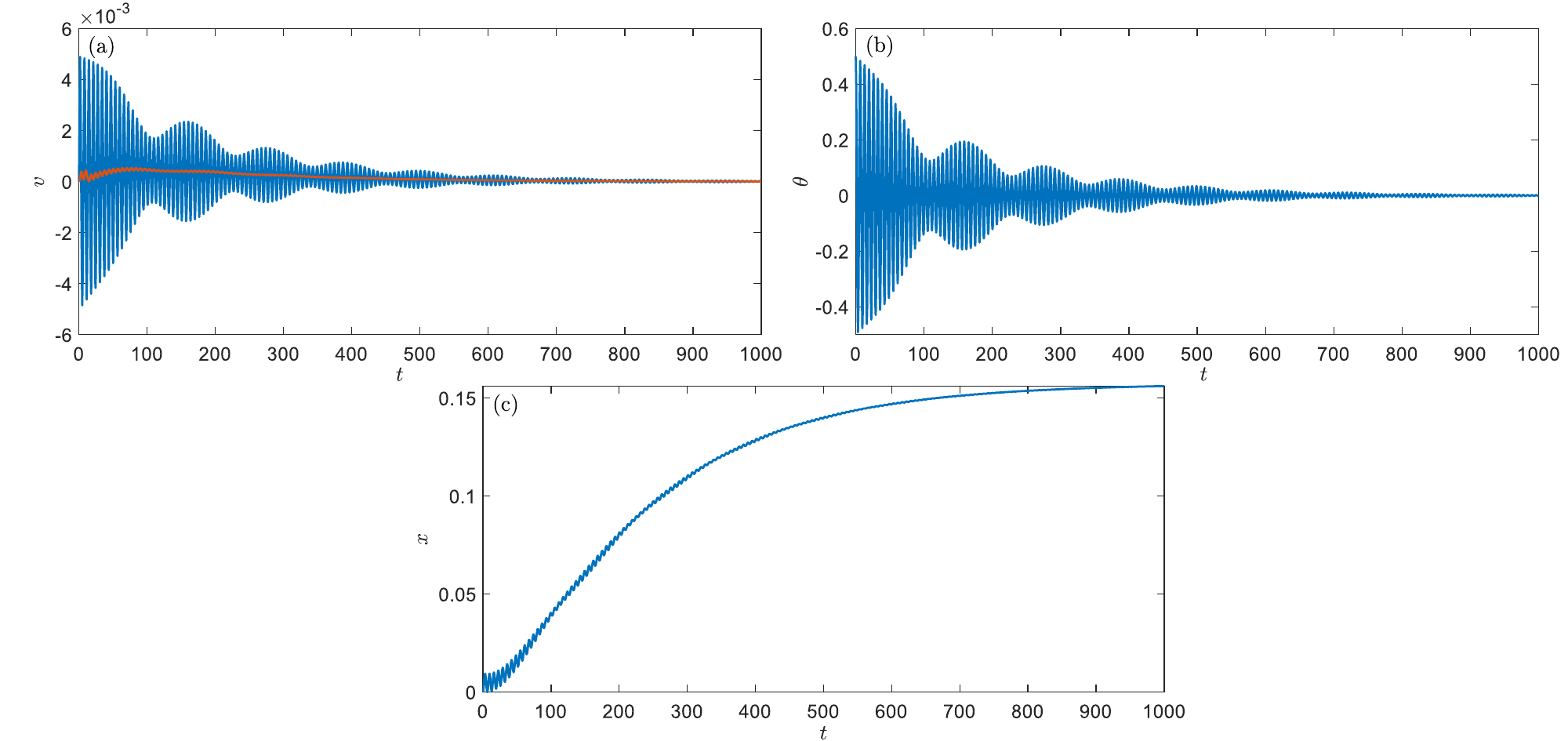}
\caption{(\textbf{Case 4}) Time histories illustrating the pendulum and capsule response under (2:1) resonance conditions, where the pendulum exhibits an oscillatory motion. System parameters: $\omega=1.94$, $\varepsilon=0.01$, $A=0.08$, $\zeta=0.01$, $\mu_1=0.01$, and $\mu_2=0.02$. Initial conditions: $x(0)=0$, $\theta(0)=0.5$, $x'(0)=0$, and $\theta'(0)=0$. (a) The time history of the capsule's velocity $v=x'$ versus time $t$, (b) the time history of the pendulum's angle $\theta$ versus time, and (c) the time histories of capsule's position $x$ versus time.}\label{F:fig13}
\end{figure*}

\begin{figure*}[ht!]
\includegraphics[width=\textwidth]{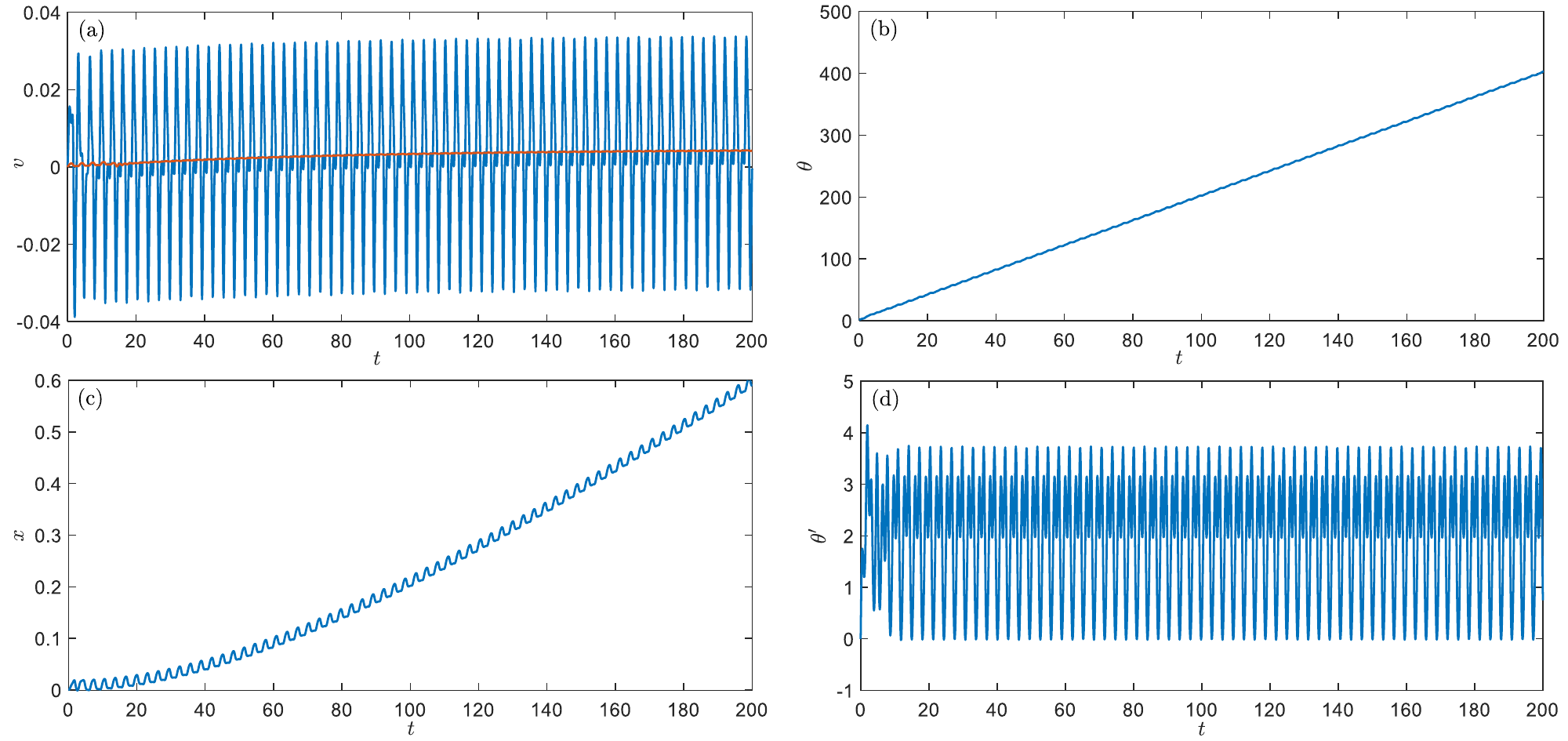}
\caption{Time histories illustrating the pendulum and capsule response under (1:1) resonance conditions, where the pendulum exhibits the perfectly rotatory response. System parameters: $\omega=2$, $\varepsilon=0.01$, $A=8$, $\zeta=1$, $\mu_1=0.01$, and $\mu_2=0.02$. Initial conditions: $x(0)=0$, $\theta(0)=2$, $x'(0)=0$, and $\theta'(0)=0$. (a) The velocity of the capsule $v=x'$ versus time $t$, (b) the angle of the pendulum $\theta$ versus time, (c) the capsule's position $x$ versus time, and (d) the temporal derivative of the pendulum angle $\theta'$ versus time. The orange curve in (a) depicts the running average of the velocity.}\label{F:fig3}
\end{figure*}

\section{Numerical evidences for the two operational modes of vibration-driven locomotion}\label{S:3}

In this section, we numerically demonstrate the existence of the resonant response regimes described above, characterized by either a purely oscillatory (2:1) response or a rotatory (1:1) response of the pendulum, accompanied by  axial progressive motion of the capsule driven by the pendulum’s internal dynamics. It is important to note that the initial angle of the pendulum in the oscillatory (2:1) regime is assumed to be small, whereas in the rotatory (1:1) regime, the pendulum is given a sufficiently large initial excitation such that resonance capture occurs near the (1:1) resonance manifold. 

In Figures~\ref{F:fig11}--\ref{F:fig13} we show time histories in a low-excitation case, namely a (2:1) oscillatory resonance of the pendulum, where in order to obtain different regimes presented in these figures, we vary the excitation frequencies or actuation amplitudes, and/or the initial angle of the pendulum. Note that each figure contains three panels: the capsule velocity $v(t)=\mathrm{d}x/\mathrm{d}t$ (with its running average), the pendulum angle $\theta(t)$, and the capsule position $x(t)$. Unless stated otherwise, the following parameter values are used throughout:
\begin{equation}\label{E:param_val}
    \zeta=0.01, \quad \mu_1=0.01, \quad \text{and} \quad \mu_2=0.02.
\end{equation}

First, in Figure~\ref{F:fig11} we present \textbf{Case~1}, using the parameter values specified in~\eqref{E:param_val} together with $A=0.01$ and $\omega=2$. The results show that after an initial transient the capsule and the pendulum stop their motion i.e., we obtain a trivial solution. Next, we consider \textbf{case 2}, with the same parameter values as in case 1, except that now the amplitude is $A=0.08$, as shown in Figure~\ref{F:fig2}. One can see that after an initial transient the capsule moves with a constant (non vanishing) average velocity (with periodic oscillations) and the pendulum oscillates about the average $\bar{\theta}=0$ with a constant amplitude. Next, we consider two additional cases, \textbf{Case~3} and \textbf{Case~4}, using the parameter values specified in~\eqref{E:param_val} together with $A=0.08$ and $\omega=1.94$. The initial conditions are $\theta(0)=2$ and $\theta(0)=0.5$, respectively. We obtain that for sufficiently large initial angle, $\theta(0)=2$ (see case 3 in Figure~\ref{F:fig12}) the regime of motion is similar to case 2, in the sense that after an initial transient the capsule moves with a constant (non vanishing) average velocity (with periodic oscillations) and the pendulum oscillates about the average $\bar{\theta}=0$ with a constant amplitude. However, for sufficiently small initial angle, $\theta(0)=0.5$ (see case 4 in Figure~\ref{F:fig13}), the system regime is similar to case 1, since after an initial transient the solution converges to the trivial one. Thus, our numerical results indicate the existence of various regimes of motion. In Section~\ref{S:41}, we use asymptotic approximations and multi-scale analysis in order to analyze the following questions analytically: what parameters affect the existence of stable solutions in various regimes, how to distinguish between various regimes (existence of bifurcation points), how many stable solutions exist within each regime, what guarantees unidirectional capsule propulsion, and how to attain maximum propulsion velocity.

The numerical results for the second, a high-excitation, case giving rise to a (1:1) rotatory resonance are shown in Figure~\ref{F:fig3}. In this Figure we present 4 panels showing the capsule velocity $v(t)$ (including the running average of the velocity), the pendulum angle $\theta(t)$, the capsule position $x(t)$, and the temporal derivative of the angle $\theta'(t)$. It can be seen that in this case, the capsule moves forward with a constant (non-vanishing) average velocity, and its angular velocity $\theta'(t)$ oscillates periodically about a nonzero mean value. The corresponding analytical study, which allows to distinguish between various regimes of motion in this case (existence of bifurcation points) is given in Section~\ref{S:42}. 

Results of Figures~\ref{F:fig2} and \ref{F:fig3} clearly show the appearances of passive progressive motion states for both aforementioned excitation cases. As it will be demonstrated below, the entire capsule dynamics based on the intrinsic mechanism of transformation of ambient structural vibrations into the progressive motion of the capsule can be fully described analytically.

In Sections~\ref{S:41} and~\ref{S:42}, we present a comprehensive analytical and numerical investigation of the intrinsic dynamics of the capsule-pendulum model governed by~\eqref{E:7}. Particular emphasis is placed on the bifurcation analysis of steady-state solutions that determine distinct locomotion regimes of the capsule, the emergence of unidirectional motion under parametric excitation, and the conditions for its existence.

\section{Analysis of oscillatory response, near (2:1) linear resonance}\label{S:41}
We now perform asymptotic analysis of the system’s oscillatory response near the (2:1) resonance. In this case we assume that
\begin{equation}\label{E:8}
\begin{aligned}
&{A} = \varepsilon {P},\quad \zeta  = \varepsilon \xi, \quad {P}  = O(1),\\
&\mu \left( {x'} \right) = \left\{ \begin{array}{l}
\varepsilon {\mu _1},\quad x' > 0,\\
\varepsilon {\mu _2},\quad x' \le 0,
\end{array} \right.  \quad {\mu_1},\,\, {\mu_2}  = O(1),
\end{aligned} 
\end{equation}
and assume the following scaling for the coordinates of pendulum angle and capsule velocity,
\begin{equation}\label{E:10}
    \theta(t)  = \sqrt \varepsilon  \Theta(t),\quad x'(t) = v(t) = \sqrt \varepsilon  V(t). 
\end{equation}

Moreover, we assume that the pendulum is subjected to a low-amplitude parametric excitation in the vicinity of the fundamental (2:1) parametric resonance. In the present case, the excitation frequency is taken to be slightly mistuned from the exact (2:1) parametric resonance condition, and is therefore expressed in the following form:

\begin{equation}\label{E:11}
    \omega  = 2 + \varepsilon \sigma ,\quad \sigma  = O\left( 1 \right).
\end{equation}

Substituting the assumptions in~\eqref{E:8}, \eqref{E:10}, and~\eqref{E:11} into the system in~\eqref{E:7}, and expanding these equations up to orders $O(\varepsilon^3)$ and $O(\varepsilon^2)$, we get that
\begin{equation}\label{E:12}
    \begin{aligned}
&{V}' = \varepsilon \left( {\Theta  - \mu(V)V} \right)\\ 
&+{\varepsilon ^2}\left({\Theta  + \Theta {{\Theta'}^2} - \frac{2}{3}{\Theta }^3 - {P}\cos \left( {\omega t } \right)\tilde \Theta  - \mu (V) V + \xi {\Theta'}}\right)\\
&+ O({\varepsilon ^3}),\\
&\Theta''  + \Theta  = \varepsilon \left( {{P}\cos \left( {\omega t } \right)\Theta  + \frac{1}{6}{{\Theta }^3} - \Theta  + \mu (V)V - \xi \Theta'} \right)\\
&+ O({\varepsilon ^2}).
\end{aligned}
\end{equation}

\subsection{Multi-scale expansion and derivation of the slow-flow model}
Following the well-known method of asymptotic complexification analysis (see e.g.~\cite{Vakakis_2008,Manevitch_2011}) we introduce the pendulum motion in the following complex form, 
\begin{equation}\label{E:13}
    {\Theta'}  + i\Theta  = \Psi.
\end{equation}
Here $\Psi$  is a complex variable which contains both the angular displacement (imaginary part) and the angular velocity (real part) of the pendulum. To analyze the vicinity of the principle parametric resonance (2:1), we exploit the method of multiple scales. Thus using the two-term multi-scale analysis we introduce the two time-scales (${t _0} = t,\,{t _1} = \varepsilon t$ ) and pursue the following expansion (see e.g.~\cite{Nayfeh_Mook_1995,Holmes_2013}),
\begin{equation}\label{E:14}
    \begin{aligned}
&\partial_t  = {\partial _{{t _0}}} + \varepsilon {\partial _{{t _1}}},\\
&V = \varepsilon {s_0}\left( {{t _0},{t _1}} \right) + {\varepsilon ^2}{s_1}\left( {{t _0},{t _1}} \right) + O\left( {{\varepsilon ^3}} \right),\\
&\Theta'  + i\Theta  = {\Psi _0}\left( {{t _0},{t _1}} \right) + \varepsilon {\Psi _1}\left( {{t _0},{t _1}} \right).
\end{aligned}
\end{equation}

Further, substituting~\eqref{E:14} into~\eqref{E:12} and collecting the terms in each order of the asymptotic expansion (i.e. $\varepsilon^0$, $\varepsilon^1$, $\varepsilon^2$), we derive the following set of equations which needs to be solved consequently starting from the lowest order $O(1)$.

\noindent
$O(1)$:
\begin{equation}\label{E:15}
    {\partial _{{t _0}}}{\Psi _0} - i{\Psi _0} = 0,
\end{equation}
\noindent
$O(\varepsilon)$:
\begin{equation}\label{E:16}
    {\partial _{{t _0}}}{s_0} =  - \frac{i}{2}\left( {{\Psi _0}({t _0},{t _1}) - {{\bar \Psi }_0}({t _0},{t _1})} \right),
\end{equation}
\begin{equation}\label{E:17}
\begin{aligned}
       &\left( {{\partial _{{t _0}}} - i} \right){\Psi _1} =  
 - {\partial _{{t _1}}}{\Psi _0}\\
 &- i\frac{{{{P}}}}{2}\left( {{e^{i\omega t_0 }} + {e^{ - i\omega t_0 }}} \right)\left( {\frac{{{{\rm{\Psi }}_0}}}{2} - \frac{{{{\bar \Psi }_0}}}{2}} \right)\\
 &- i\left( {\frac{{{{\rm{\Psi }}_0}}}{2} - \frac{{{{\bar \Psi }_0}}}{2}} \right)-\xi \left( {\frac{{{{\rm{\Psi }}_0}}}{2} + \frac{{{{\bar \Psi }_0}}}{2}} \right)-\\
 &- \frac{i}{{48}}\left( {{{\bar \Psi }_0}^3 - 3{{\bar \Psi }_0}{{\left| {{{\rm{\Psi }}_0}} \right|}^2} + 3{{\rm{\Psi }}_0}{{\left| {{{\rm{\Psi }}_0}} \right|}^2} - {\rm{\Psi }}_0^3} \right), 
 \end{aligned}
\end{equation}
$O(\varepsilon^2)$:
\begin{equation}\label{E:18}
\begin{aligned}
    &{\partial _{{t _1}}}{s_0} + {\partial _{{t _0}}}{s_1} = 
 - i\left( {\frac{{{{\rm{\Psi }}_1}}}{2} - \frac{{{{{\rm{\bar \Psi }}}_1}}}{2}} \right) - i\left( {\frac{{{{\rm{\Psi }}_0}}}{2} - \frac{{{{{\rm{\bar \Psi }}}_0}}}{2}} \right)\\
 &+ \frac{i}{8}\left( { - {\rm{\Psi }}_0^3 - {{\rm{\Psi }}_0}{{\left| {{{\rm{\Psi }}_0}} \right|}^2} + {{{\rm{\bar \Psi }}}_0}{{\left| {{{\rm{\Psi }}_0}} \right|}^2} + {{{\rm{\bar \Psi }}}_0}^3} \right)\\
 &- \frac{i}{{12}}\left( {3{{{\rm{\bar \Psi }}}_0}{{\left| {{{\rm{\Psi }}_0}} \right|}^2} - 3{{\rm{\Psi }}_0}{{\left| {{{\rm{\Psi }}_0}} \right|}^2} - {{{\rm{\bar \Psi }}}_0}^3 + {\rm{\Psi }}_0^3} \right)\\
 &+ i\frac{P}{2}\left( {{e^{i\omega t_0 }} + {e^{ - i\omega t_0} }} \right)\left( {\frac{{{{\rm{\Psi }}_0}}}{2} - \frac{{{{{\rm{\bar \Psi }}}_0}}}{2}} \right)\\
 &+ \xi\left( {\frac{{{{\rm{\Psi }}_0}}}{2} + \frac{{{{{\rm{\bar \Psi }}}_0}}}{2}} \right) - {\rm{\mu }}\left( {{s_0}} \right){s_0}.
 \end{aligned}
\end{equation}

Solving equation~\eqref{E:15} with respect to the fast time scale, we derive the general solution for the pendulum motion in the following complex form,
\begin{equation}\label{E:19}
    {\Psi _0} = \varphi \left( {{t _1}} \right){e^{i{t _0}}},
\end{equation}
where it is possible to express $\varphi(t_1)$ in polar form, 
\begin{equation}\label{E:20}
    \varphi(t_1)  = a\left( {{t _1}} \right){e^{i\delta \left( {{t _1}} \right)}}.
\end{equation}
Proceeding to the next order of approximation, inserting~\eqref{E:19} into~\eqref{E:16} and using the compact trigonometric representation, we arrive at the following equation
\begin{equation}\label{E:21}
    {\partial _{{t _0}}}{s_0} = \left| {\varphi \left( {{t _1}} \right)} \right|\sin \left( {{t _0} + \delta \left( {{t _1}} \right)} \right).
\end{equation}
Integrating~\eqref{E:21} once, with respect to the fast time scale ($t_0$), we derive the general form of the capsule response, reading
\begin{equation}\label{E:22}
    s_0(t_0,t_1) = D\left( {{t _1}} \right) - \left| {\varphi \left( {{t _1}} \right)} \right|\cos \left( {{t _0} + \delta \left( {{t _1}} \right)} \right).
\end{equation}
Plugging~\eqref{E:11} and~\eqref{E:19} into~\eqref{E:17}, we arrive at:
\begin{equation}\label{E:23}
    \begin{aligned}
    &{\partial _{{t _0}}}{\Psi _1} - i{\Psi _1} =  - (\partial _{{t _0}}\varphi){e^{i{t _0}}} + 
\Bigl[\frac{i}{2}\left( {\varphi {e^{i{t _0}}} - \bar \varphi {e^{ - i{t _0}}}} \right)\\
&- \frac{i{P}}{4}\left( {{e^{i\left( {2{t _0} + \sigma {t _1}} \right)}} + {e^{ - i\left( {2{t _0} + \sigma {t _1}} \right)}}} \right)\left( {\varphi {e^{i{t _0}}} - \bar \varphi {e^{ - i{t _0}}}} \right)\\
 &+ \frac{1}{6}{\left( { - \frac{i}{2}\left( {\varphi {e^{i{t _0}}} - \bar \varphi {e^{ - i{t _0}}}} \right)} \right)^3} -  \frac{\xi}{2}\left( {\varphi {e^{i{t _0}}} + \bar \varphi {e^{ - i{t _0}}}} \right)\Bigr].
 \end{aligned}
\end{equation}
The sum of all resonant (secular) terms multiplying $e^{i t_0}$ must vanish in order to enforce the solution to be bounded and thus to satisfy the asymptotic scaling assumption. Hence, collecting the secular terms in the RHS of~\eqref{E:23} and equating to zero, we arrive at the following complex ODE depicting the complex amplitude evolution of the pendulum response,
\begin{equation}\label{E:24}
    \partial_{t_1}\varphi = \left( {\frac{i}{2}\varphi  + \frac{i}{4}{P}\bar \varphi {e^{i\sigma {t _1}}} - \frac{i}{{16}}{{\left| \varphi  \right|}^2}\varphi  - \frac{{\xi }}{2}\varphi } \right). 
\end{equation}
To bring~\eqref{E:24} into autonomous form, we introduce additional transformation
\begin{equation}\label{E:25}
    \varphi  = \phi {e^{\frac{{i\sigma {t _1}}}{2}}}
\end{equation}
and rewrite~\eqref{E:24} in the following form,
\begin{equation}\label{E:26}
    {\phi _{{t _1}}} = \frac{i}{2}\left( {1 - \sigma } \right)\phi  + \frac{i}{4}{P}\bar\phi  - \frac{i}{{16}}{\left| \phi  \right|^2}\phi  - \frac{\xi}{2}\phi.
\end{equation}

To derive the analytical approximation for the capsule motion, we  proceed to the final order of asymptotic expansion~\eqref{E:18}. Using the Fourier expansion in~\eqref{E:22}, we Fourier-expand the dissipative term with respect to the DC term and the first fast-time scale harmonic, yielding
\begin{equation}\label{E:27}
    \mu(V)V = {d_0}\left( {{t_1}} \right) + {d_2}\left( {{t _1}} \right)\cos \left( {{t _0} + \delta \left( {{t _1}} \right)} \right),
\end{equation}
where $d_0$ and $d_2$ depend on $\mu_1$, $\mu_2$, $|\varphi|$, and $D$.
All the details of the Fourier expansion including the functional form of the Fourier coefficients $d_0$, $d_1$ can be found in Appendix A. Substituting~\eqref{E:27} into~\eqref{E:18}, we have that
\begin{equation}\label{E:28}
\begin{aligned}
    {\partial _{{t _0}}}{s_1} &=-\Bigl[ {\partial _{{t _1}}}D + {d_0}\Bigr]  - 
\left| \phi  \right|{\partial _{{t _1}}}\delta \sin \left( {{t _0} + \delta } \right)\\
&- \left( {{d_2} - {\partial _{{t _1}}}\left| \phi  \right|} \right)\cos \left( {{t _0} + \delta } \right)\\
  &- i\left( {\frac{{{{\rm{\Psi }}_1}}}{2} - \frac{{{{{\rm{\bar \Psi }}}_1}}}{2}} \right) - i\left( {\frac{{{{\rm{\Psi }}_0}}}{2} - \frac{{{{{\rm{\bar \Psi }}}_0}}}{2}} \right)\\
 &+ \frac{i}{8}\left( { - {\rm{\Psi }}_0^3 - {{\rm{\Psi }}_0}{{\left| {{{\rm{\Psi }}_0}} \right|}^2} + {{{\rm{\bar \Psi }}}_0}{{\left| {{{\rm{\Psi }}_0}} \right|}^2} + {{{\rm{\bar \Psi }}}_0}^3} \right)\\
 &- \frac{i}{{12}}\left( {3{{{\rm{\bar \Psi }}}_0}{{\left| {{{\rm{\Psi }}_0}} \right|}^2} - 3{{\rm{\Psi }}_0}{{\left| {{{\rm{\Psi }}_0}} \right|}^2} - {{{\rm{\bar \Psi }}}_0}^3 + {\rm{\Psi }}_0^3} \right)\\
 &+ i\frac{P}{2}\left( {{e^{i\left( {2{t _0} +  \sigma t_1 } \right)}} + {e^{ - i\left( {2{t _0} +  \sigma t_1} \right)}}} \right)\left( {\frac{{{{\rm{\Psi }}_0}}}{2} - \frac{{{{{\rm{\bar \Psi }}}_0}}}{2}} \right)\\
 &+ \xi\left( {\frac{{{{\rm{\Psi }}_0}}}{2} + \frac{{{{{\rm{\bar \Psi }}}_0}}}{2}} \right).
\end{aligned}
\end{equation}

Note that the ``secular'' term $[\partial_{t_1} D + d_0]$ in~\eqref{E:28} creates unbounded terms diverging in time, in contradiction to the scaling assumption. In order to avoid this divergence, these terms are required to vanish. Hence, by eliminating the secular terms in the RHS of~\eqref{E:28}, we derive the following slow-evolution equation:
\begin{equation}\label{E:29}
    {\partial _{{t _1}}}D =  - {d_0}\left( {D,\left| \varphi  \right|} \right),
\end{equation}
where the derivation of the explicit expression for $d_0$ is given in Appendix A. Equation~\eqref{E:29} determines the slow evolution of the DC component (i.e. the slow drift) of the capsule velocity. This component represents the average of capsule velocity (see e.g. Figure~\ref{F:fig2}(a)). 
We close the present subsection with summarizing the modulation equations depicting the slow-amplitude dynamics of oscillatory response of the pendulum as well as the slow drift of capsule velocity,
\begin{equation}\label{E:30}
    \begin{aligned}
&{\partial _{{t _1}}}\phi  = \frac{i}{2}\left( {1 - \sigma } \right)\phi  + \frac{i}{4}{P}\bar \phi  - \frac{i}{{16}}{\left| \phi  \right|^2}\phi  - \frac{\xi}{2}\phi, \\
&{\partial _{{t _1}}}D =  - {d_0}\left( {D,\left| \phi  \right|} \right).
\end{aligned}
\end{equation}

\subsection{Analysis of stationary-points of the slow-flow}
We start the analytical treatment of the slow-flow given in \eqref{E:30} from the analysis of its stationary points. Thus, seeking for the fixed points of the reduced slow-flow model we set all the time derivatives to be zero (${\partial _{{t _1}}}D = {\partial _{{t _1}}}\phi  = 0$). This leads to the following set of nonlinear equations:
\begin{subequations}\label{E:31}
   \begin{align}
&\frac{i}{2}\left( {1 - \sigma } \right)\phi  + \frac{i}{4}{P}\bar\phi  - \frac{i}{{16}}{\left| \phi  \right|^2}\phi  - \frac{\xi}{2}\phi  = 0, \label{E:31A}\\
&{\left( {{{\rm{\mu }}_1} - {{\rm{\mu }}_2}} \right)}\sqrt {{{\left| \phi  \right|}^2} - {D^2}} \label{E:31B}\\
&\quad\quad+ {D}\left( {\left( {{{\rm{\mu }}_2} - {{\rm{\mu }}_1}} \right)\arccos \left( {\frac{{D}}{{\left| \phi  \right|}}} \right) + {\rm{\pi }}{{\rm{\mu }}_1}} \right) = 0. \notag
\end{align}
\end{subequations}

Performing standard algebraic manipulations, we solve \eqref{E:31A} and obtain the following three distinct branches of solutions,
\pagebreak
\begin{equation}\label{E:32}
    \begin{aligned}
&{\left| \phi_0  \right|} = 0,\\
&{\left| \phi_1  \right|} = \sqrt {8\left( {1 - \sigma } \right) + 4\sqrt {P^2 - {{\left( {2\xi } \right)}^2}} }\\
&\quad\quad{\rm{ whenever}} \quad {P} > 2\xi,\quad \sigma \leq 1 + \frac{1}{2}\sqrt {P^2 - 4{{\xi}^2}}, \\
&{\left| \phi_2  \right|} = \sqrt {8\left( {1 - \sigma } \right) - 4\sqrt {P^2 - {{\left( {2\xi } \right)}^2}} }\\
&\quad\quad{\rm{ whenever}} \quad {P} > 2\xi, \quad \sigma  \leq 1 - \frac{1}{2}\sqrt {P^2 - 4{{\xi }^2}}.
\end{aligned}
\end{equation}
For further details regarding the solution of~\eqref{E:31A}, see Appendix B. 

Substituting the distinct branches of the steady-state solutions $\left| \phi_0  \right|$, $\left| \phi_1  \right|$, and $\left| \phi_2  \right|$ corresponding to  the pendulum motion into the second equation of~\eqref{E:31} and solving for $D$, we derive implicitly the three branches $D_0$, $D_1$, and $D_2$ of the steady state solutions corresponding the averaged capsule velocities, respectively,
\begin{equation}\label{E:33}
    \begin{aligned}
&{{D}_0} = 0,\\
&{\left( {{{\rm{\mu }}_1} - {{\rm{\mu }}_2}} \right)}\sqrt {\left| \phi_{1,2}  \right|^2 - {{D}_{1,2}}}\\
&\phantom{{{\rm{\mu }}_2}}\,\,\,\,+ {{D}_{1,2}}\left( {\left( {{{\rm{\mu }}_2} - {{\rm{\mu }}_1}} \right)\arccos \left( {\frac{{{{D}_{1,2}}}}{{{{\left| \phi_{1,2}  \right|}}}}} \right) + {\rm{\pi }}{{\rm{\mu }}_1}} \right) = 0.
\end{aligned}
\end{equation}
Note that although the equations for $D_{1}(|\phi_1|)$ and $D_{2}(|\phi_2|)$ are transcendental, in the relevant region both equations have a unique solution.

\begin{figure*}[ht!]
\includegraphics[width=\textwidth]{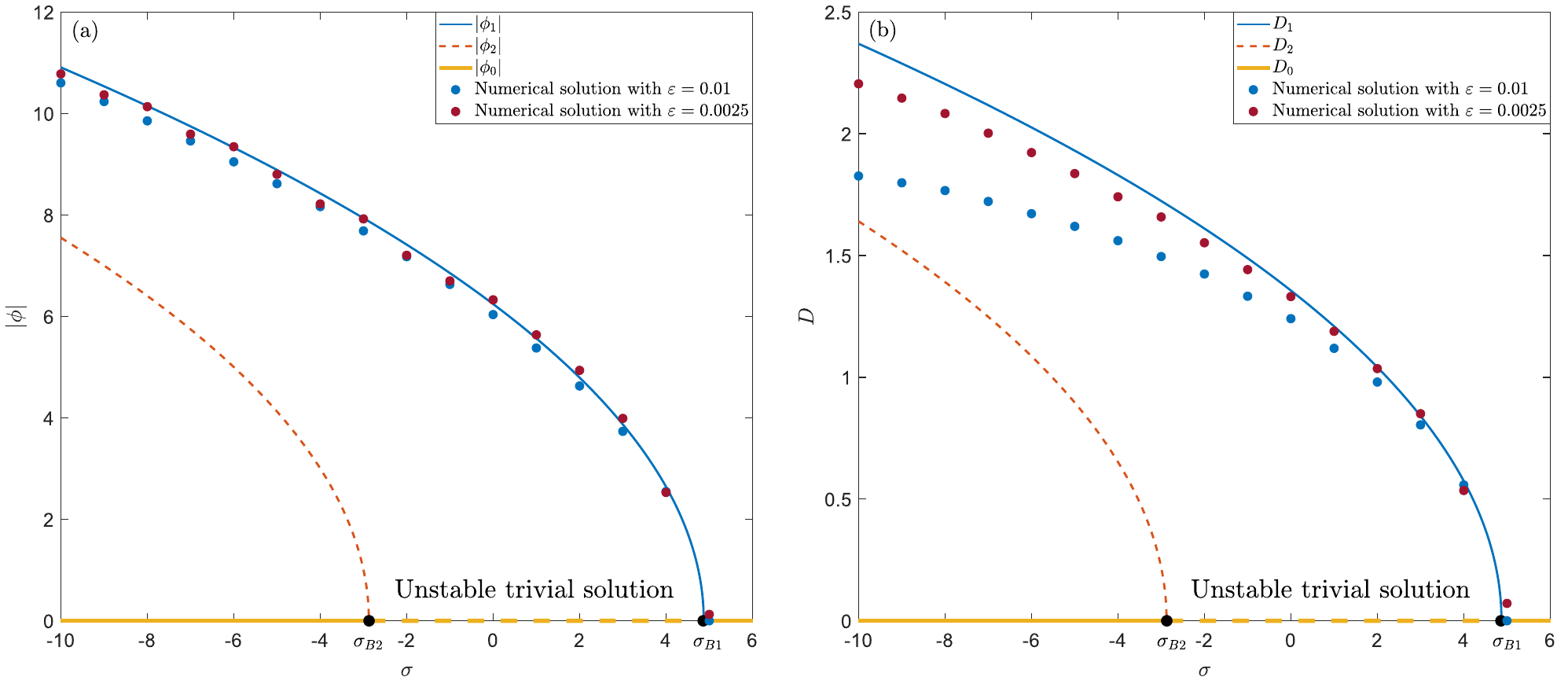}
\caption{Bifurcation diagrams for the oscillatory response with the parameters: $P=8$, $\xi=1$, $\mu_1=1$, and $\mu_2=2$ for (a) pendulum amplitudes $|\phi_0|$, $|\phi_1|$, and $|\phi_2|$ and (b) the corresponding capsule's average velocities $D_0$, $D_1$, and $D_2$. Solid curves represent stable steady state regimes, dashed curves represent the unstable ones, while the blue and the maroon dots represent the numerical results obtained by solving the system of equation in~\eqref{E:7} with $\varepsilon=0.01$  and $\varepsilon=0.0025$, respectively.}\label{F:fig4}
\end{figure*}

Bifurcation diagrams corresponding to different solution branches of the amplitude $|\phi|$ of the steady-state response of internal pendulum as well as the steady-state average velocities $D$ of the capsule model are illustrated in Figure~\ref{F:fig4}.

\begin{figure*}[ht!]
\includegraphics[width=\textwidth]{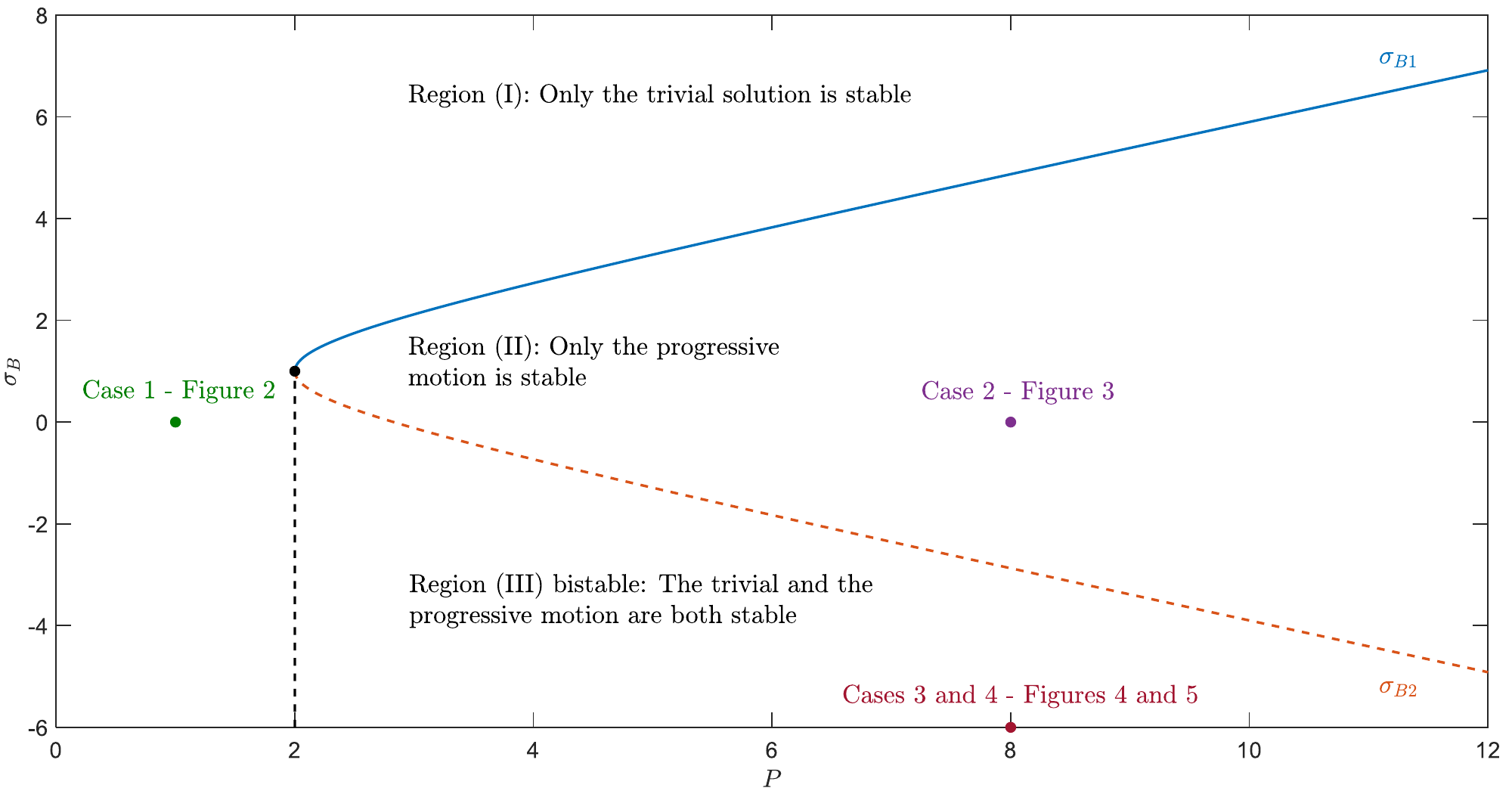}
\caption{Stability transition curves in $(P,\sigma)$ plane for the oscillatory response with the parameter values $\xi=1$, $\mu_1=1$, and $\mu_2=2$. The full circles visualize the simulated cases, which were shown in Section~\ref{S:3}.}\label{F:fig14}
\end{figure*}

\begin{figure*}[ht!]
\includegraphics[width=\textwidth]{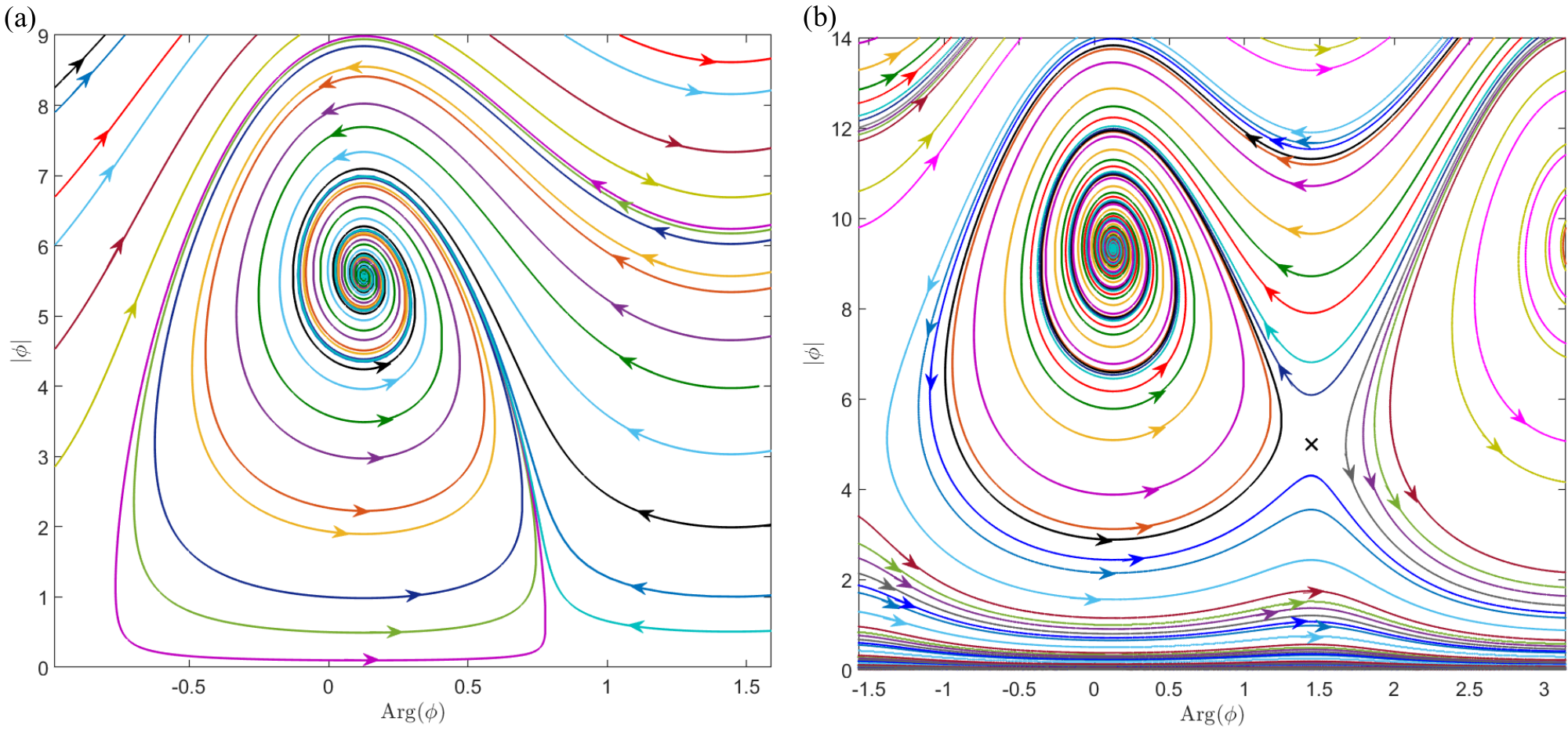}
\caption{Phase planes of the slow-flow corresponding to the oscillatory response. System parameters: $P=8$, $\xi=1$, $\mu_1=1$, and $\mu_2=2$ for (a) $\sigma_{B2}<\sigma=1<\sigma_{B1}$ (b) $\sigma=-6<\sigma_{B2}$. The arrows indicate the direction of propagation (as time increases) of the various solutions (each solution marked by its own color) and `x' in panel (b) denotes the saddle point}.\label{F:fig5}
\end{figure*}

The $\sigma$ axis of the diagrams corresponds to the frequency detuning in the vicinity of the fundamental parametric resonance (2:1). The solid lines indicate the stable steady state regimes while dashed lines indicate the unstable ones. The blue and maroon dots represent the numerical results obtained by solving the system in~\eqref{E:7} with $\varepsilon=0.01$ and $\varepsilon=0.0025$, respectively. It can be seen that for larger values of $\sigma$, the analytical solution agrees well with the numerical one, but as the values of $\sigma$ decrease, the discrepancy increases, especially in panel (b). Moreover, on can observe that decreasing $\varepsilon$, the difference  between the analytical approximation and the numerical decreases as well, so that the error is approximately proportional to $\varepsilon$. As is clear from the bifurcation diagrams the system undergoes two subsequent bifurcations at the points 
\[{\rm{ }}{\sigma _{B1}} = 1 + \frac{1}{2}\sqrt{P^2 - 4{\xi }^2}, \quad \left( {P > 2\xi} \right)\] 
and 
\[{\sigma _{B2}} = 1 - \frac{1}{2}\sqrt {P^2 - 4{{\xi }^2}}, \quad \left( {P > 2\xi} \right),\]
which split the system dynamics into the three qualitatively different regions. 

These three regions in $(P,\sigma)$ plane denoted by Region (I), Region (II), and Region (III), as well as the stability transition curves are shown in Figure~\ref{F:fig14}. For sufficiently large $\sigma$ or sufficiently small $P$, we are in Region (I), where only trivial solutions are stable (see case 1-Figure~\ref{F:fig11}). Within the region bounded by the solid blue and dashed orange curves, namely Region (II), only the solution corresponding to progressive capsule motion is stable, as illustrated by a representative example in case 2-Figure~\ref{F:fig2}. The region (III) surrounded by dashed orange and black curves is bistable in the sense that both solutions, the trivial and progressive capsule motion solutions, are stable, and we can obtain either solution based on the initial conditions, as illustrated in cases 3 and 4 (see Figures~\ref{F:fig12}-\ref{F:fig13}).

In Figure~\ref{F:fig5}, we plot $|\phi|$ versus $\mathrm{Arg}(\phi)$ obtained by solving equation~\eqref{E:26}; hereafter, graphs of this type are referred to as ``phase diagrams''. The figure illustrates the attraction of trajectories to various steady-state solutions whose structure differs qualitatively across the parameter regions. Thus, in the first region $\left(\sigma > \sigma_{B1}\right)$, which is not shown in the figure, there exists a single stable trivial solution, manifested by the inevitable decay of the system to the zero state for any initial condition.

Obviously, this state of dynamical system is highly undesirable as the desired stationary locomotion state cannot be achieved in this region. As the detuning parameter decreases and passes through  $\left( {\sigma  = {\sigma _{B1}}} \right)$ the first super-critical pitchfork bifurcation occurs. In this case the trivial solution (zero state) becomes unstable, and a pair of new stable fixed points emerge. These stable fixed points of the slow-flow model correspond to a stable, time periodic (steady state) solution of the original model. As is evident from the bifurcation diagram, by choosing the frequency detuning parameter $\sigma$ within the interval $\left( \sigma_{B2} < \sigma < \sigma_{B1} \right)$, one can readily conclude that, for any arbitrary initial condition, the system response converges to the nontrivial stationary state associated with the desired progressive motion. The corresponding slow-flow dynamics of the pendulum reduced to the plane, is shown in Figure~\ref{F:fig5}(a). Additional third region, which deserves a separate consideration corresponds to the range of detuning parameter $\left( {\sigma< {\sigma _{B2}}} \right)$. As the detuning parameter decreases and passes through the second sub-critical pitchfork bifurcation point the unstable trivial fixed point switches to a stable focus, and an additional pair of unstable fixed points (saddle-points) emerge (the saddle point is marked by `x').
As is clear from the phase plane of Figure~\ref{F:fig5}(b), for low initial amplitude $|\phi|$, the solution will inevitably decay to a zero state, while, for sufficiently large initial values of $|\phi|$, the trajectories will converge to a non-trivial fixed point of the phase plane $(\text{Arg}(\phi),|\phi|)$, which reflects stable steady state oscillations of parametric pendulum. 

\begin{figure*}[ht!]
\includegraphics[width=\textwidth]{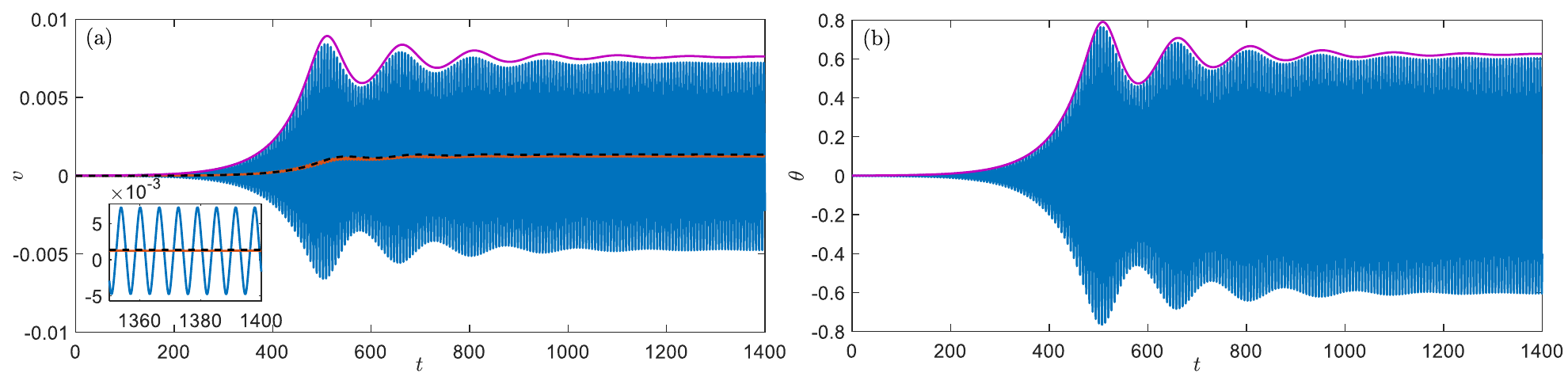}
\caption{Comparison of the original model (blue solid line) and the envelope obtained by the analytical approximation given in~\eqref{E:AV} and~\eqref{E:ATHETA} (violet) for (a) the capsule velocity and (b) the pendulum angle, respectively, for the following parameter values: $\omega=2$, $\varepsilon=0.01$, $A=8$, $\zeta=1$, $\mu_1=1$, and $\mu_2=2$, with the initial conditions: $x(t=0)=0$, $\theta(t=0)=0.001$, $x'(t=0)=0$, and $\theta'(t=0)=0$, whereas for the slow-flow model we used $\phi(t_1=0)=(0.001/\sqrt{\varepsilon})i$ and $D(t_1=0)=0$. The orange curve in (a) denotes the running average of the numerical solution and the dashed black curve represents the mean velocity obtained by integrating the slow-flow model in~\eqref{E:mean_V}. In the inset, we show an increased view of the graph.}\label{F:fig9}
\end{figure*}

From the perspective of practical system design, we aim in this study to derive the fastest possible progressive motion of the capsule. As is clear from the bifurcation analysis brought above, by choosing the system's parameters in the Region II the system response will inevitably converge to the stationary locomotion state. This region, indeed guarantees us the emergence of the stationary forward motion of the capsule for any choice of initial conditions. By inspecting the average velocity of the capsule, we can see from the diagrams that their values become higher as the frequency detuning parameter $\sigma$ decreases. However, unlike the Region II, Region III is a bit more involved. On one hand, the higher velocity values of the capsule locomotion obtained in this region make this region advantageous in comparison to Regions I and II. On the other hand, the emergence of the stable trivial solution still retains the possibility for the termination of the permanent progressive motion of the capsule under low initial excitations of internal device.

\subsection{Numerical verification}\label{S:Num ver first part}
To confirm the validity of the developed analytical model corresponding to the oscillatory response of the pendulum we compare the results (i.e. time - series) of numerical simulations of the original model~\eqref{E:7} with these of the slow-flow model~\eqref{E:30}. 

The amplitude envelope denoted with the violet solid line on both panels of Figure~\ref{F:fig9} corresponds to the modulated amplitudes of the capsule velocity and pendulum oscillations, which are given by 
\begin{align}
&\text{A}_V(t)=\varepsilon^{3/2}[D(t_1/\varepsilon)+\phi(t_1/\varepsilon)], \label{E:AV}\\
&\text{A}_{\theta}(t)=\sqrt{\varepsilon}\phi(t_1/\varepsilon), \label{E:ATHETA}
\end{align}
respectively, (see Figure~\ref{F:fig9}), where $\phi$ and $D$ were found according to the stable branches of~\eqref{E:32} and~\eqref{E:33}. Blue solid line denotes the response of the original model in equation~\eqref{E:7}, orange curve is the running average of the numerical results, and the dashed black line is the average of the capsule's velocity calculated by
\begin{equation}\label{E:mean_V}
    \bar{V}(t)=\varepsilon^{3/2}D(t_1/\varepsilon).
\end{equation}
 As is clear from these comparison plots, the agreement between the true response of the original model and the modulation equations is fairly satisfactory.

\section{Analysis of rotatory response, nonlinear (1:1) resonance}\label{S:42}
This section is devoted to asymptotic analysis of the rotatory response demonstrated numerically in Section~\ref{S:3} (see Figure~\ref{F:fig3}). In contrast to the first type of locomotion regime, where only a weak initial perturbation of the pendulum was assumed, the strongly nonlinear locomotion regime considered here is initiated by a large initial excitation of the internal pendulum. This ensures that the system begins in the vicinity of the (1:1) nonlinear resonance manifold, which means that the actuation frequency is not necessarily equal to the natural frequency of linear oscillations, and undergoes immediate resonance capture.

For the rotatory (1:1) resonance, we adopt the following asymptotic scaling
\begin{equation}\label{E:9}
    \begin{aligned}
&{A} = O(1),\quad {\zeta}  = O(1),\\
&\mu \left( {x'} \right) = \left\{ \begin{array}{l}
\varepsilon {\mu _1},\quad x' > 0,\\
\varepsilon {\mu _2},\quad x' \le 0.
\end{array} \right.
\end{aligned}
\end{equation}

\subsection{Averaged-flow model near (1:1) parametric resonance}\label{S:421}
In this section we apply direct averaging method and thus derive the averaged-flow model in the vicinity of (1:1) parametric resonance.
Substituting the asymptotic assumptions and asymptotic scaling made in~\eqref{E:9}  into~\eqref{E:7}, we obtain
\begin{equation}\label{E:34}
    \begin{aligned}
&v' = \varepsilon {{\theta'}^2}\sin (\theta ) - \varepsilon \cos (\theta )\theta''  - \varepsilon \mu (v)v,\\
&\theta''  =  - {\zeta}\theta'  - \left( {1 - A\cos (\omega t )} \right)\sin (\theta ) - v'\cos (\theta ).
\end{aligned}
\end{equation}

\begin{figure*}[ht!]
\includegraphics[width=\textwidth]{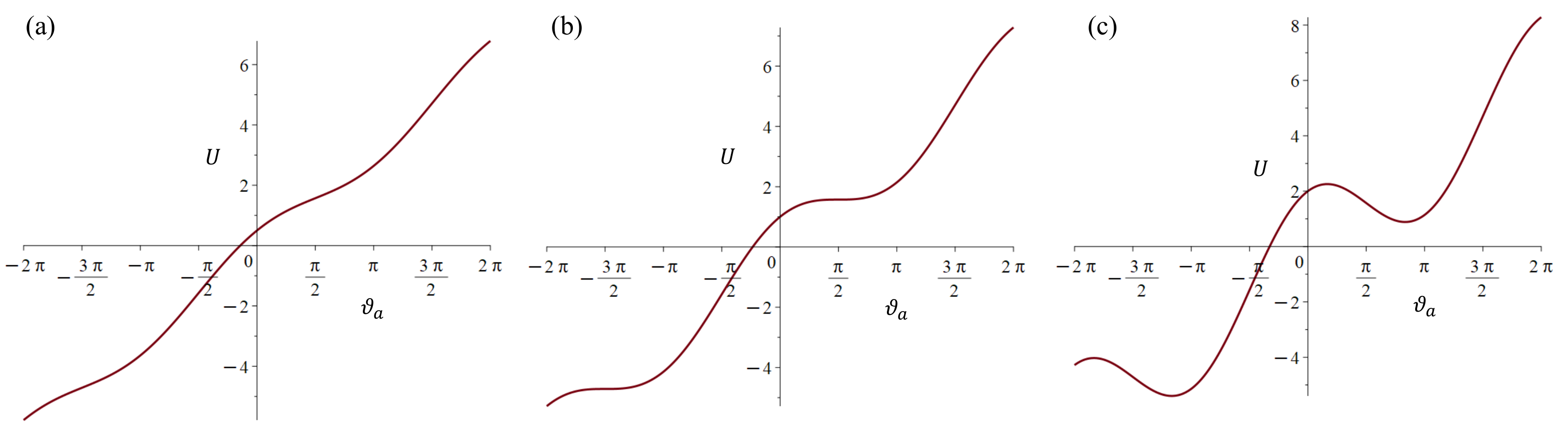}
\caption{Graphs of $U$ versus $\vartheta_a$ for (a) $\eta=0.5$, (b) $\eta=1$, and (c) $\eta=2$.}\label{F:fig6}
\end{figure*}

Plugging the first equation of~\eqref{E:34} into the second one, we arrive at the following set of equations,
\begin{equation}\label{E:35}
    \begin{aligned}
&v' = \varepsilon {{\theta'}^2}\sin (\theta ) - \varepsilon \cos (\theta )\theta''  - \varepsilon \mu (v)v,\\
&\theta''  + \left( {1 - A\cos (\omega t )} \right)\sin (\theta ) =  - {\zeta} \theta'\\
&\phantom{\ddot \theta}+ \varepsilon \left( {{{\cos }^2}(\theta )\theta''  + \mu (v)v\cos (\theta ) - \frac{1}{2}{{\theta'}^2}\sin (2\theta )} \right).
\end{aligned}
\end{equation}

We formulate the following ansatz for the motion of internal pendulum,
\begin{equation}\label{E:36}
    \theta  = \omega t  + \vartheta (t ),
\end{equation}
where $\vartheta(t)$ is the slow varying resonance phase.
Substituting~\eqref{E:36} into~\eqref{E:35}, neglecting the $O(\varepsilon)$ terms and averaging over a single forcing period $T=2\pi/\omega$, one obtains the following equation of motion describing the drift of the averaged phase $\vartheta$ of the rotating pendulum,
\begin{equation}\label{E:37}
    \vartheta''_a  =  - {\zeta} \left( {\omega  + \vartheta'_a} \right) + \frac{A}{2}\sin \left( {{\vartheta_a}} \right),
\end{equation}
where $\vartheta_a=\langle \vartheta(t)\rangle=(1/T)\int_{0}^{T}\vartheta(t)\,dt$ denotes the phase of the pendulum angle averaged over one period of the forcing. The integration is carried out over a single excitation period, during which the fast oscillatory component of $\vartheta(t)$ is averaged out, while the slowly varying component $\vartheta_a(t)$ (being treated as a constant in the integration process) is retained after the integration.

We look for the conditions on ${\vartheta_a}$ as a function of ${A}$, ${\zeta}$, and $\omega$ that coincide with the assumption of the rotatory resonant motion of the pendulum (i.e. $\theta  = \omega t  + {\vartheta}(t )$).
Note that equation~\eqref{E:37} can be represented as follows,
\begin{equation}\label{E:38}
    \vartheta''_a  + {\xi} \vartheta'_a  =  - \frac{d U}{d\vartheta_a},
\end{equation}
where the derivative $dU/d\vartheta_a$ is given by,
\begin{equation}\label{E:39}
    \frac{d U}{d\vartheta_a} = \zeta \omega  - \frac{{{A}}}{2}\sin \left( {\vartheta_a} \right).
\end{equation}
Integrating~\eqref{E:39} with respect to $\vartheta$, we derive the potential function,
\begin{equation}\label{E:40}
    U\left( {\vartheta_a} \right) = {\zeta} \omega \vartheta_a  + \frac{{{{A}}}}{2}\cos \left( {\vartheta_a} \right).
\end{equation}

Now, we find the conditions for the convergence of $\vartheta_a$ to a constant value, which is ``phase locking'' of the average solution with the excitation frequency. From~\eqref{E:39}, we conclude that in order to enforce this to occur, the potential $U$ should have a minimum point, which corresponds to a locally stable equilibrium of average phase $\vartheta_a$. However, $U$ is a transcendental function in $\vartheta_a$ that has a monotonic and an oscillatory parts, and thus the condition for local minima existence is determined by the ratio between these two parts, $\eta=A/(2{\zeta}\omega)$, where
\begin{equation}\label{E:41}
    \frac{{U\left( { \vartheta_a } \right)}}{{{\zeta} \omega }} = \vartheta_a  + \eta \cos \left( {\vartheta_a} \right).
\end{equation}
It can be shown that for $\eta>1$ the potential $U(\vartheta_a)$ is monotonic, and for $\eta<1$, $U(\vartheta_a)$ has two minimum points. The critical transition value $\eta=1$ undergoes a saddle-node bifurcation between the state with minimum and the state with monotonic increase.

In Figure~\ref{F:fig6}, we illustrate the three possible forms of the potential function $U(\vartheta_a)$, where panels (a)-(c) correspond three following cases:
\begin{itemize}
    \item [(a)] $U(\vartheta_a)$ is monotonically increasing. This implies the absence of phase locking which leaves no possibility for the sustained (1:1) resonance.
    \item [(b)] Point of the saddle-node bifurcation.
    \item [(c)] Existence of two minimum points for $U(\vartheta_a)$, which guarantees the existence of two solutions of phase locking, leading to the (1:1) resonant rotatory response.
\end{itemize}

\begin{figure}[ht!]
\includegraphics[width=0.48\textwidth]{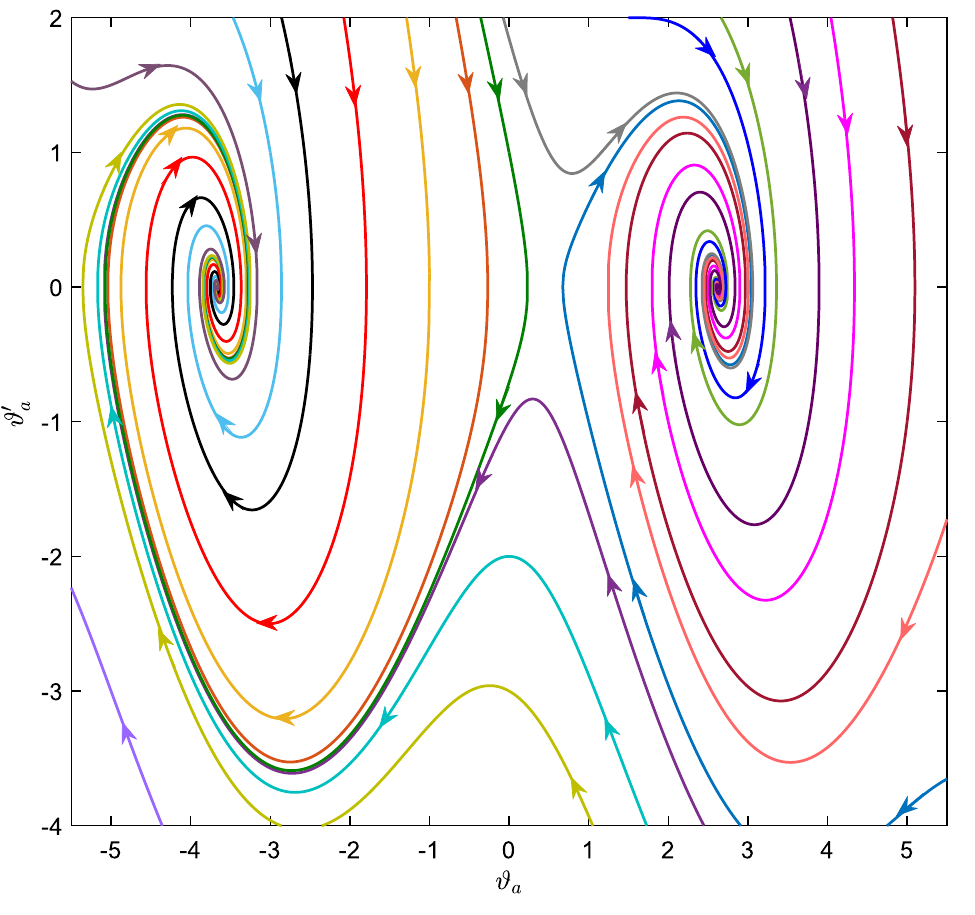}
\caption{Phase plane of the averaged-flow corresponding to the (1:1) rotatory response. System parameters: $A=8$, $\zeta=1$, and $\omega=2$ (so that $\eta>1$). The arrows indicate the direction of propagation (as time increases) of the various solutions (each solution marked by its own color).}\label{F:fig7}
\end{figure}

In Figure~\ref{F:fig7}, we show the corresponding phase plane for the case of existence of the phase locking domain which brings the possibility for (1:1) resonance capture. In this figure we focus on $\eta>1$. 

The trajectories in Figure~\ref{F:fig7} are obtained by numerically integrating equation~\eqref{E:37} for various initial conditions. As shown, all trajectories converge to one of the two phase-locking points of the (1:1) nonlinear resonance (i.e. this resonance corresponds to the rotational motion of the pendulum), which correspond to the minima of $U(\vartheta_a)$ in Figure~\ref{F:fig6}. Moreover, this phase diagram corresponds to the numerical simulation presented in Figure~\ref{F:fig3}, which depicts the time histories of the solution in the case of (1:1) rotational resonance motion of the pendulum.

Thus, picking the system parameters such that $\eta>1$ and choosing the initial conditions inside the phase locking domain we guarantee the emergence of (1:1) resonant locomotion states manifested by the rotatory motion of the internal pendulum. 

To establish the averaged velocities of the capsule being engaged in this resonant regime, we proceed with the direct averaging procedure. This procedure assumes that the capsule velocity, can be separated into the oscillatory like part $\tilde{u}$ (i.e. the component of the response which does not include the drift in capsule velocity) as well as the non-oscillatory component $D$ containing this drift (i.e. the averaged part of the velocity). Using this form we formulate the two-component representation of the capsule velocity,
\begin{equation}\label{E:42}
    v = D + \tilde u,
\end{equation}
where we request that $\tilde u$ is an oscillatory function with zero mean value.
Plugging~\eqref{E:42} along with the ansatz adopted for the resonant motion of the pendulum~\eqref{E:36} into the first equation of~\eqref{E:34}, we arrive at the following equation,
\begin{equation}\label{E:43}
\begin{aligned}
    D' + \tilde{u}' &= \varepsilon {\left( {\omega  + \vartheta'_a } \right)^2}\sin (\omega t  +  \vartheta_a )\\
    &- \varepsilon \cos (\omega t  + \vartheta_a)\vartheta''_a  - \varepsilon \mu\left((D + \tilde u)\right)\left( {D + \tilde u} \right).
    \end{aligned}
\end{equation}
The equation~\eqref{E:43} contains the non-oscillatory component as well as the oscillatory-like part $\tilde{u}$. Next, we split the equation  into the drifting component part, yielding
\begin{equation}\label{E:44}
    D' =  - \varepsilon \mu\left((D + \tilde u)\right)\left( {D + \tilde u} \right)
\end{equation}
as well as the oscillatory part, which yields the following equation,
\begin{equation}\label{E:45}
    \begin{aligned}
    \tilde{u}' &= \varepsilon \Bigl[ {{\left( {\omega  + {\vartheta'_a}} \right)}^2}\sin (\omega t  +  \vartheta_a )\\
    &- \left( {\zeta} \left( {\omega  + {\vartheta'_a}} \right) - \frac{{{{A}}}}{2}\sin \left( {\vartheta_a } \right) \right)\cos (\omega t  + \vartheta_a ) \Bigr].
    \end{aligned}
\end{equation}
Note that in the derivation of~\eqref{E:45} we used~\eqref{E:37}. To find $\tilde{u}$, we integrate~\eqref{E:45} with respect to $t$, under the assumption of slow drift of the averaged phase terms (i.e. terms being treated as constants in the integration process) as well as zero initial conditions ($\tilde{u}(0)=0$),
\begin{equation}\label{E:46}
    \begin{aligned}
    \tilde u(t) &=  - \frac{\varepsilon }{\omega }\Bigl[ {{\left( {\omega  + \vartheta'_a } \right)}^2}\cos (\omega t  + \vartheta_a )\\
    &- \left( {\zeta \left( {\omega  + \vartheta'_a } \right) - \frac{A}{2}\sin \left( {\vartheta_a} \right)} \right)\sin (\omega t  + \vartheta_a) \Bigr].
    \end{aligned}
\end{equation}
This solution can be represented in the following, compact trigonometric form,
\begin{equation}\label{E:47}
    \tilde u =  - \varepsilon B\left( {\vartheta_a,\vartheta'_a } \right)\cos \left( {\omega t  + \phi_a \left( {\vartheta_a,\vartheta'_a} \right)} \right),
\end{equation}
where the amplitude $B( {\vartheta_a ,{\vartheta'_a} })$ and phase coefficients \linebreak ${\phi_a}(\vartheta_a,\vartheta'_a)$ are given by
\begin{equation}\label{E:48}
    \begin{aligned}
    B\left( {\vartheta_a ,\vartheta'_a} \right) &= \frac{1}{\omega }\Bigl[{{\left( {\omega  + \vartheta'_a}\right)}^4}\\
    &+ {{\left( {{\xi} \left( {\omega  + \vartheta'_a } \right) - \frac{{A}}{2}\sin \left( {\vartheta_a} \right)} \right)}^2}\Bigr]^{1/2}
    \end{aligned}
\end{equation}
\begin{equation}\label{E:49}
    \phi_a \left( {\vartheta_a ,\vartheta'_a} \right) = \vartheta_a + \arctan \left[ {\frac{{{\zeta}\left( {\omega  + \vartheta'_a} \right) - \frac{{{A}}}{2}\sin \left( {\vartheta_a} \right)}}{{{{\left( {\omega  + \vartheta'_a} \right)}^2}}}} \right].
\end{equation}


Returning to the second component $D$ of ~\eqref{E:44}, we substitute expression~\eqref{E:47} into~\eqref{E:44} and apply the direct averaging procedure. This yields the following slow evolution equation for the capsule velocity–drift component:
\begin{equation}\label{E:51}
    \begin{aligned}
    D'_a  =  - \frac{{\varepsilon \omega }}{{2\pi }}\int_0^{\frac{{2\pi }}{\omega }}&\Bigl[\mu \left( {{D_a} - \varepsilon B\cos \left( {\omega t  + \phi_a } \right)} \right)\\
    &\times
    \left( {{D_a}  - \varepsilon B\cos \left( {\omega t  + \phi_a } \right)} \right)\Bigr] dt.
    \end{aligned}
\end{equation}
For convenience, we define the new phase as
\begin{equation}\label{E:52}
    \psi  = \omega t  + \phi.
\end{equation}
We can see that over one periodic cycle the damping coefficient has jumps according to the following relation, reformulated in terms of two-component representation,
\begin{equation}\label{E:53}
     \mu \left( {{D_a} + \tilde u} \right) = \left\{ \begin{array}{l}
{{\mu}_1},  \quad   {D_a} + \tilde u > 0,\\
{{\mu}_2},  \quad   {D_a} + \tilde u \le 0.
\end{array} \right.
\end{equation}
It can be easily shown that over a single oscillatory period and under the assumption of ${D_a} < \varepsilon B(\vartheta_a ,\vartheta'_a)$, the function \linebreak ${D_a} - \varepsilon B(\vartheta_a ,\vartheta'_a)\cos \psi$ switches between the different values of damping parameter at three points of $\psi$:
\begin{equation}\label{E:54}
    \begin{aligned}
    &{\psi _1} =  - \arccos\left( {\frac{{D_a}}{{\varepsilon B\left( {\vartheta_a ,\vartheta'_a} \right)}}} \right), \\
    &{\psi _2} = \arccos\left( {\frac{D_a}{{\varepsilon B\left( {\vartheta_a ,\vartheta'_a} \right)}}} \right), \\
    &{\psi _3} = 2\pi  - \arccos\left( {\frac{D}{{\varepsilon B\left( {\vartheta_a ,\vartheta'_a } \right)}}} \right).
    \end{aligned}
\end{equation}
It is worthwhile noting the following restriction:\linebreak $- \frac{\pi }{2} \le \arccos \left( {\frac{D_a}{{\varepsilon B\left( {\vartheta_a ,\vartheta'_a} \right)}}} \right) \le \frac{\pi }{2}$.  Substituting~\eqref{E:52} and~\eqref{E:54} into~\eqref{E:51} yields the following averaged part (averaged capsule velocity) of the $D_a$ component,
\begin{equation}\label{E:55}
    \begin{aligned}
    D'_a   =  &- \frac{\varepsilon }{{2\pi }}\Bigl( \int_{{\psi _1}}^{{\psi _2}} {\mu_2}\left( D_a - \varepsilon B\cos \psi  \right)d\psi\\
    &+ \int_{{\psi _2}}^{{\psi _3}} \mu_1\left(D_a - \varepsilon B\cos \psi  \right)d\psi  \Bigr).
    \end{aligned}
\end{equation}
Hence, the evolution equation for the averaged velocity, obtained from the integration of~\eqref{E:55}, is given by
\begin{equation}\label{E:56}
    \begin{aligned}
    D'_a =  &- \frac{{\rm{\varepsilon }}}{{\rm{\pi }}}\Bigl[ \left( {{{{\rm{\mu }}}_1} - {{{\rm{\mu }}}_2}} \right)\left( {D_a\arcsin \left( {\frac{D_a}{{{\rm{\varepsilon }}B}}} \right) + \sqrt {{{\left( {{\rm{\varepsilon }}B} \right)}^2} - {D_a}^2} } \right)\\
    &+ \frac{{{\rm{\pi }}\left( {{{{\rm{\mu }}}_1} + {{{\rm{\mu }}}_2}} \right)}}{2} D_a \Bigr].
    \end{aligned}
\end{equation}

The choice of the initial conditions satisfies the relation between the pendulum and the phase drift, respectively, 
\begin{equation}\label{E:57}
    \begin{aligned}
&\theta (0) =   \vartheta_a (0),\\
&\theta'(0) = \omega + \vartheta'_a(0).
\end{aligned}
\end{equation}
According to the problem statement, the initial capsule velocity is assumed to be zero, i.e. 
\begin{equation}\label{E:58}
    v(0) = 0.
\end{equation}
The initial condition in~\eqref{E:58} can be represented by the oscillatory and D.C. components 
\begin{equation}\label{E:59}
    D(0) + \tilde u(0) = 0.
\end{equation}
From ~\eqref{E:59} we directly obtain the relation between the initial amplitude and phase of the oscillating component of the capsule's velocity to the initial value of its drifting component. This relation, reads 
\begin{equation}\label{E:60}
    \begin{aligned}
&D(0) =  - \tilde u(0) =  - \left( {0 - \varepsilon B(0)\cos \left( {\phi_a \left( 0 \right)} \right)} \right),\\
&D(0) = \varepsilon B(0)\cos \left( { \phi_a \left( 0 \right)} \right).
\end{aligned}
\end{equation}

\subsection{Analysis of stationary points of the averaged-flow}
To find the amplitude of the steady state regimes of the pendulum - capsule model we set the time derivatives, in equations~\eqref{E:37} and~\eqref{E:56}, to zero and obtain the following set of algebraic equations
\begin{equation}\label{E:61}
    {\zeta} \omega  - \frac{A}{2}\sin \left( {\vartheta_a} \right) = 0,
\end{equation}
\begin{equation}\label{E:62}
    \begin{aligned}
     & \left( {{{{\rm{\mu }}}_1} - {{{\rm{\mu }}}_2}} \right)\left( {D_a}\arcsin \left( {\frac{D_a} {{{\rm{\varepsilon }}B}}} \right) + \sqrt {{{\left( {{\rm{\varepsilon }}B} \right)}^2} - {D_a}^2}  \right)\\
     &+ \frac{{D_a}\pi\left(\mu_1 + \mu _2 \right)}{2} = 0.
     \end{aligned}
\end{equation}

The stationary points of $\vartheta_a$, are easily obtained from~\eqref{E:61},
\begin{equation}\label{E:63}
    \begin{aligned}
    &{\vartheta _{ss}} = \arcsin\left( {\frac{{2{\zeta} \omega }}{{{A}}}} \right) +  n\pi,\quad n = 0,\pm{2},\pm{4}\ldots,\\
    &{\vartheta _{ss}} =  n\pi-\arcsin\left( {\frac{{2{\zeta}\omega }}{{{A}}}} \right),\quad n = \pm{1},\pm{3}\ldots.
    \end{aligned}
\end{equation}
Note that the stationary points of $\vartheta_a$, are not the equilibrium points of the pendulum angle, but rather they are the equilibrium points of the average phase difference between pendulum oscillations and periodic excitation.
Further, using ~\eqref{E:63} in~\eqref{E:62}, we obtain that the stationary points of ${D_a}$ are given explicitly by the following nonlinear transcendental algebraic equation,
\begin{equation}\label{E:62A}
    \begin{aligned}
     & \left( {{{{\rm{\mu }}}_1} - {{{\rm{\mu }}}_2}} \right)\left( {D_{ss}}\arcsin \left( {\frac{D_{ss}} {{{\rm{\varepsilon }}B({\vartheta_{ss}})}}} \right) + \sqrt {{{\left( {{\rm{\varepsilon }}B({\vartheta_{ss}})} \right)}^2} - {D_{ss}}^2}  \right)\\
     &+ \frac{D_{ss}\pi\left(\mu_1 + \mu _2 \right)}{2} = 0.
     \end{aligned}
\end{equation}

\begin{figure*}[ht!]
\includegraphics[width=\textwidth]{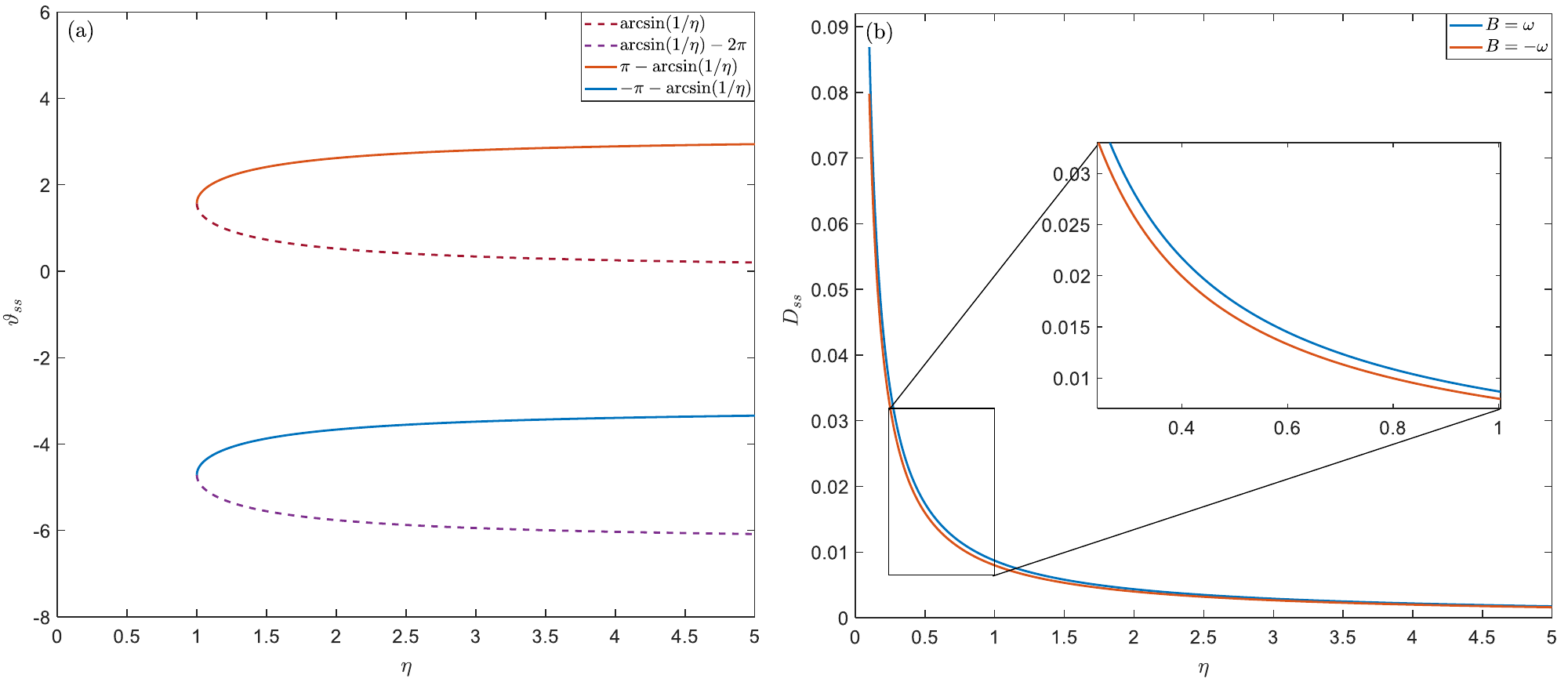}
\caption{Bifurcation diagrams corresponding to the (1:1) rotatory response. System parameters: $A=8$, $\zeta=1$, $\mu_1=1$,  $\mu_2=2$, and $\omega=A/(2\zeta\eta)$. (a) Rotatory phase drift response $\vartheta_{ss}$ and (b) averaged capsule velocity $D_{ss}$ versus $\eta$. The solid curves in (a) represent stable solutions and dashed curves represent unstable solutions. The inset in (b) shows the magnified view of the region indicated with a rectangle.}\label{F:fig8}
\end{figure*}

Note that substituting~\eqref{E:61} into~\eqref{E:48} and assuming that $\vartheta'_a$ vanishes, we get that the stationary values of $B$ corresponding to the same fixed point solutions of the rotatory locomotion states are given by
\begin{equation}\label{E:64}
    B\left( {{{\vartheta }_{ss}}} \right) =\pm\omega.
\end{equation}

The stationary points in~\eqref{E:63} are stable for odd values of $n$, since for odd values of $n$ the potential $U(\vartheta_a)$ attains minimum at $\vartheta_a=\vartheta_{ss}$, while for even values of $n$ the potential $U(\vartheta_a)$ attains maximum at $\vartheta_a=\vartheta_{ss}$. In Figure~\ref{F:fig8}, we show the bifurcation diagram for the steady-state solutions of the resonant phase $\vartheta _{ss}$ and the steady state solutions of the averaged capsule velocity ${D}_{ss}$ as functions of the system parameter $\eta={A}/(2{\zeta}\omega)$ calculated by using equations~\eqref{E:63} and~\eqref{E:62A} with $B(\vartheta_{ss})$ substituted from~\eqref{E:64}. 

\begin{figure*}[ht!]
\includegraphics[width=\textwidth]{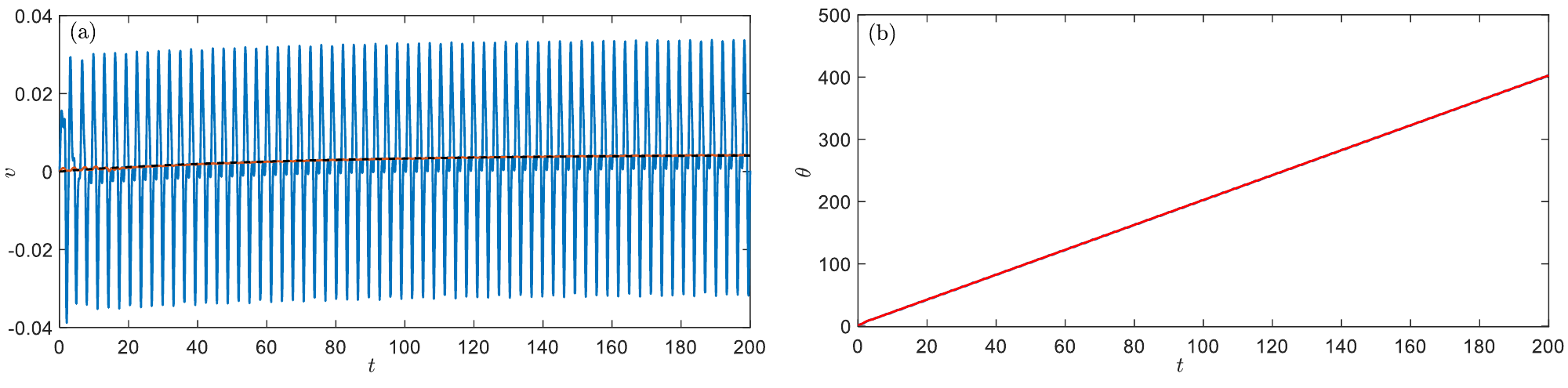}
\caption{Comparison of numerical solution of the original model (blue solid line) and the running average obtained by the analytical approximation  for (a) the capsule velocity and (b) the pendulum angle for the following parameter values: $\omega=2$, $\varepsilon=0.01$, $A=8$, $\zeta=1$, $\mu_1=1$, and $\mu_2=2$, with the initial conditions: $x(0)=0$, $\theta(0)=2$, $x'(0)=0$, and $\theta'(0)=0$. The orange curve in (a) denotes the running average of the numerical solution and the dashed black curve represents the mean velocity obtained by integrating the averaged-flow equation, given in~\eqref{E:56} for $D_a$ with $D_a(0)=0$. The red curve in (b) denotes the mean angle obtained by using~\eqref{E:36} with $\vartheta_a$ that was found by integrating the averaged-flow equation in~\eqref{E:37} subject to the initial conditions $\vartheta_a(0)=2$, $\vartheta'_a(0)=-2$.}\label{F:fig10}
\end{figure*}

As is clear from the results of Figure~\ref{F:fig8}(a), a pair of nontrivial locomotion states emerges through a saddle-node bifurcation which occurs exactly at $\eta = 1$. Phase planes illustrating the typical averaged-flow system dynamics of the averaged pendulum phase below and above the threshold are brought in Figure~\ref{F:fig7}(a) and Figure~\ref{F:fig7}(b), respectively. As expected below the threshold (Figure~\ref{F:fig7}(a)) the system converges to zero while above the threshold (Figure~\ref{F:fig7}(b)) the system evolves towards a stable focus corresponding to a stable, non-trivial locomotion state. An additional noteworthy result is the distinct difference in the capsule’s average propulsion speed when the pendulum rotates clockwise versus counterclockwise. This directional asymmetry is clearly reflected in the bifurcation diagram presented in Figure~\ref{F:fig8}(b).

\subsection{Numerical verification}
In Figure~\ref{F:fig10}, we show again the numerical results that were presented in Figure~\ref{F:fig3}, where now we compare them with the theoretical results derived in this section. Specifically, for the case of the capsule velocity response, we compare the moving average of the capsule velocity derived from the time histories of the original model~\eqref{E:35} with that of the averaged $D_a$ component of the averaged-flow~\eqref{E:56} (see Figure~\ref{F:fig10}(a)). As for the pendulum motion, similarly to the case of oscillatory response of internal pendulum illustrated in Section~\ref{S:Num ver first part}, we compare the averaged amplitude of the pendulum response given by the averaged-flow model~\eqref{E:37} with the time-series plot of the true response of the original model~\eqref{E:35}. To guarantee the occurrence of a resonance capture leading to the permanent (1:1) resonant rotatory response we choose the proper initial excitation for the internal pendulum which reads ($\theta \left( 0 \right) = 0,\, \theta' \left( 0 \right) = \omega $) along with the proper choice of the system parameters $\eta  = A/(2\zeta\omega) > 1$  (in accordance with the analysis in  Section~\ref{S:421}). 

As is clear from the time histories plots of Figure~\ref{F:fig10} for the case of the rotatory response, we observe a rather satisfactory agreement between the results of the analytically derived averaged model with the numerical computation of moving average out of the original signal. It is also clear that the special choice of initial conditions for the pendulum is in accordance with the phase drift of the pendulum, meaning that for the special choice of initial excitation applied on the internal pendulum the system will exhibit the internal resonant rotatory response resulting in the progressive motion of the capsule model whose velocity is dictated by the three system parameters, namely, the amplitude of excitation  ($A$), internal dissipation level ($\zeta$), and excitation frequency ($\omega$).

\section{Conclusions}\label{S:6}
In the present work, we examined the dynamics of a two-degree-of-freedom capsule model composed of a capsule incorporating an internal pendulum. The capsule is supported by a vertically oscillating foundation and is free to move horizontally, subject to asymmetric viscous damping. A dedicated analytical treatment was developed for the system response: a multi-scale analysis for the principal parametric (2:1) resonance, and a direct averaging procedure for the rotatory (1:1) resonance.
For the oscillatory regime, the analysis of the resulting slow-flow model and its bifurcation diagram reveals the existence of progressive motion state under two distinct nonlinear response scenarios. Depending on the detuning parameter, the capsule may exhibit:
(i) a robust, uniquely attracting steady progressive motion state, characterized by a constant positive average velocity independent of initial conditions; or
(ii) a higher-velocity progressive motion state, whose realization depends on the initial displacement of the pendulum and corresponds to an additional nontrivial stationary solution emerging through a saddle-node bifurcation.
For the rotatory regime, the averaged averaged-flow model and its bifurcation structure indicate a qualitatively different progressive motion mechanism. Here, sustained rotation of the pendulum accompanied by forward capsule motion is achievable only when the parameter $\eta$ is greater than unity, and when the initial conditions place the system within the domain of attraction of the corresponding stable stationary point (see Figure~\ref{F:fig7}). In this regime, the averaged capsule velocity is determined solely by the excitation frequency and damping parameters and attains a strictly positive constant value.
The analytical predictions show very good agreement with the numerical simulations.
Several directions for further research naturally emerge from this study. In particular, the present analysis focused on a piecewise-linear damping law. Extending the methodology to more complex nonlinear dissipation mechanisms exhibiting stick-slip transitions would be of considerable interest and may lead to new realizations or enhanced implementations of pendulum-driven directional progressive motion systems. 


\section*{Declarations}

\noindent
{\bf{Funding}}
The work of Y. O. has been
supported by Israel Science Foundation under grant no. 1382/23. Y. S. acknowledges the support from the Israel Science Foundation under grant number 999/25.

\bigskip

\noindent{\bf{Conflict of interests}}
The authors declare that they have no conflicts of interest.

\bigskip

\noindent{\bf{Data availability}}
The datasets generated and/or analysed during the current study are available from the corresponding author on reasonable request. 

\bigskip

\noindent{\bf{AI-assisted writing disclosure}}
Y. S. employed AI-assisted language tools for grammar and phrasing in the Introduction, Abstract, and portions of the main text, as well as for reference identification. The authors reviewed and validated all content and assume full responsibility for the manuscript.

\bigskip

\noindent{\bf{Authors' contributions}}
Y. O. and Y. S. - conceptualization. G. I. and A. Z. - analysis and compilation of the significant part of the paper. Y. S. - compilation of the Abstract and Introduction of the manuscript. Y. O. and Y. S. edited the manuscript, supervised the research in all of its stages and acquired the funding.

\section*{Appendix A - Fourier expansion of $\mu(V)V$  with the use of the multiple scales expansion}
\setcounter{equation}{0}
\renewcommand{\theequation}{A.\arabic{equation}}
Let us consider the phase 
\begin{equation}\label{E:66}
    \psi  = {t _0} + \delta.
\end{equation}
The non-linear viscous damping of $\mu(V)V$  can be expanded up through order $O(\varepsilon^2)$ in the following manner, where $s_0$ is the solution given in~\eqref{E:22}
\begin{equation}\label{E:65}
    \begin{aligned}
&\mu \left( {{s_0}} \right){s_0} = {d_0}\left( {\left| \varphi  \right|,D} \right) + {d_1}\left( {\left| \varphi  \right|,D} \right)\sin \psi\\
&\phantom{\mu \left( {{s_0}} \right){s_0}}+ {d_2}\left( {\left| \varphi  \right|,D} \right)\cos \psi,\\
&\mu \left( s_0 \right) = \left\{ \begin{array}{l}
{\mu_1},\quad {s_0} > 0,\\
{\mu_2},\quad {s_0} \le 0.
\end{array} \right.
\end{aligned}
\end{equation}
It can be easily shown that over a single oscillatory period, the function in~\eqref{E:22} switches between the different values of damping parameter at three points of  $\psi$:
\begin{equation}\label{E:67}
    \begin{aligned}
    &{\psi _1} =  - \arccos \left( {\frac{D}{{\left| \varphi  \right|}}} \right),\\
    &{\psi _2} = \arccos \left( {\frac{D}{{\left| \varphi  \right|}}} \right),\\
    &{\psi _3} = 2\pi  - \arccos \left( {\frac{D}{{\left| \varphi  \right|}}} \right).
    \end{aligned}
\end{equation}

We restrict ourselves here to the range of \linebreak $-\pi \le \arccos \left( {\frac{D}{{\left| \varphi  \right|}}} \right) \le \pi $. Thus, the explicit representation of the coefficients for the harmonic expansion in~\eqref{E:65} can be acquired using the Fourier series coefficient integral over one periodic cycle,
\begin{equation}\label{E:68}
    \begin{aligned}
&{d_0}\left( {\left| \varphi  \right|,D} \right) = \frac{1}{{\rm{\pi }}}\Bigl[\left( \mu_1 - \mu_2 \right)\sqrt {{{\left| {\rm{\varphi }} \right|}^2} - {D^2}}\\
&\phantom{{d_0}\Bigl( {\left| \varphi  \right|,D} \Bigr)}+ D\Bigl( {\left( {{{\rm{\mu }}_2} - {{\rm{\mu }}_1}} \right)\arccos \Bigl( {\frac{{D}}{{\left| {\rm{\varphi }} \right|}}} \Bigr) + {\rm{\pi }}{{\rm{\mu }}_1}} \Bigr) \Bigr],\\
&{d_1}\left( {\left| \varphi  \right|,D} \right) = 0,\\
&{d_2}\left( {\left| \varphi  \right|,D} \right) = \frac{1}{{\left| {\rm{\varphi }} \right|{\rm{\pi }}}}\Bigl[{D}\left( {{{\rm{\mu }}_2} - {{\rm{\mu }}_1}} \right)\sqrt {{{\left| {\rm{\varphi }} \right|}^2}{\rm{ - }}{{D}^2}}\\
&\phantom{{d_2}\left( {\left| \varphi  \right|,D} \right)}+ {{\left| {\rm{\varphi }} \right|}^2}\Bigl( {\bigl( 2\mu_1 - \mu_2 \bigr)\arccos \Bigl( {\frac{{D}}{{\left| {\rm{\varphi }} \right|}}} \Bigr)} \Bigr) - 2{\rm{\pi }}|\varphi|^2{{\rm{\mu }}_1}\Bigr].
\end{aligned}
\end{equation}


We will not take into account any further harmonics; this is due to the fact that we will be deriving the slow varying components by use of secular terms which are regarded to principle parametric resonance $\omega$. The Fourier expansion in symbolic form can be expressed as,
\begin{equation}\label{E:69}
    \mu (V)V = {d_0} + {d_2}\cos \psi.
\end{equation}

\section*{Appendix B - Analytical solution of the slow-flow given in~\eqref{E:31A}}

\setcounter{equation}{0}
\renewcommand{\theequation}{B.\arabic{equation}}

In order to solve the equation in~\eqref{E:31A}, we rearrange it in the following form
\begin{equation}\label{E:70}
   \begin{aligned}
&8i\left( {1 - \sigma } \right)\phi   - i{\left| \phi  \right|^2}\phi  - 8\xi\phi  = -4iP\bar{\phi},
\end{aligned}
\end{equation}
and write the conjugate of this equation, namely
\begin{equation}\label{E:71}
   \begin{aligned}
&-8i\left( {1 - \sigma } \right)\bar{\phi}   + i{\left| \phi  \right|^2}\bar{\phi}  - 8\xi\bar{\phi}  = 4i{P}\phi.
\end{aligned}
\end{equation}
Multiplying equation~\eqref{E:70} by~\eqref{E:71}, we arrive at
\begin{equation}
    \begin{aligned}
        |\phi|^6&-16(1-\sigma)|\phi|^4+64\left[(1-\sigma)^2+\xi^2\right]|\phi|^2\\
        &=16P^2|\phi|^2.
    \end{aligned}
\end{equation}
Solving this equation, we obtain the solutions which appear in~\eqref{E:32}. 


\end{document}